\documentclass[times,sort&compress,3p]{elsarticle}
\journal{Journal of Multivariate Analysis}
\usepackage[labelfont=bf]{caption}

\usepackage{amsmath,amsfonts,amssymb,amsthm,booktabs,color,epsfig,graphicx,hyperref,url}

\usepackage{mathrsfs}
\RequirePackage[OT1]{fontenc}

\allowdisplaybreaks 

\numberwithin{equation}{section}

\theoremstyle{plain}
\newtheorem{theorem}{Theorem}
\newtheorem{lemma}{Lemma}

\newtheorem{proposition}{Proposition}

\theoremstyle{definition} 

\newtheorem{remark}{Remark}

\newcommand{\ignore}[1]{} 

\allowdisplaybreaks

\numberwithin{equation}{section} 

\font\msbmx=msbm10                   
\font\msbmvii=msbm7                  
\font\msbmv=msbm5
\def\varnothing{\mathchoice{\mbox{\msbmx\char'077}}%
{\mbox{\msbmx\char'077}}{\mbox{\msbmvii\char'077}}{\mbox{\msbmv\char'077}}}%

\def\vect{\mathop{\rm Vect}\limits}
\def\tr{\mathop{\rm tr}\limits}

\newcommand{\vae}{\varepsilon}

\newcommand{\wh}{\widehat}
\newcommand{\wt}{\widetilde}

\def\bbn{{\mathbb N}}

\def\bbr{{\mathbb R}}

\newcommand{\mb}[1]{\mathbf{#1}} 

\newcommand{\Xb}{{\mb{X}}}

\newcommand{\mc}[1]{\mathcal{#1}} 
\newcommand{\Rc}{{\mc{R}}}

\newcommand{\Nc}{{\mc{N}}}
\newcommand{\Fc}{{\mc{F}}}

\newcommand{\Hc}{{\mc{H}}}

\newcommand{\cI}{{\mc{I}}}

\newcommand{\cM}{{\mc{M}}}
\newcommand{\cF}{{\mc{F}}}
\newcommand{\cS}{{\mc{S}}}

\newcommand{\Xc}{{\mc{X}}}
\newcommand{\Bc}{{\mc{B}}}

\newcommand{\ao}{\mbox{o}}

\def\text#1{\hbox{#1}}
\def\proof{{\noindent \bf Proof. }}

\def\endproof{\mbox{\ $\qed$}}

\def\E{{\bf E}}
\def\cB{{\cal B}}

\def\K{{\bf K}}
\def\P{{\bf P}}
\def\C{{\bf C}}
\def\L{{\bf L}}
\def\p{{\bf p}}
\def\D{{\bf D}}
\def\B{{\bf B}}

\def\A{{\bf A}}
\def\U{{\bf U}}
\def\H{{\bf H}}
\def\V{{\bf V}}

\def\Q{{\bf Q}}

\def\c{{\bf c}}

\def\a{{\bf a}}
\def\e{{\bf e}}
\def\b{{\bf b}}

\def\k{{\bf k}}
\def\n{{\bf n}}
\def\g{{\bf g}}

\def\l{{\bf l}}
\def\v{{\bf v}}

\def\m{{\bf m}}
\def\q{{\bf q}}
\def\Chi{{\bf 1}}
\def\d{\mathrm{d}}
\def\build #1_#2{\mathrel{\mathop{\kern 0pt #1}\limits_\zs{#2}}}
\newcommand{\zs}[1]{{\mathchoice{#1}{#1}{\lower.25ex\hbox{$\scriptstyle#1$}}
{\lower0.25ex\hbox{$\scriptscriptstyle#1$}}}}

\numberwithin{equation}{section}
\def\proof{{\noindent \bf Proof. }}
\def\endproof{\mbox{\ $\qed$}}

\usepackage{mathrsfs}
\newcommand{\Pb}{{\mathsf{P}}} 
\newcommand{\EV}{{\mathsf{E}}} 
\newcommand{\PFA}{\mathsf{PFA}}
\newcommand{\PMI}{\mathsf{PMI}}

\newcommand{\mrm}[1]{\mathrm{#1}}

\newcommand{\F}{{\mrm{F}}}

\def\One{\mathchoice{\rm 1\mskip-4.2mu l}{\rm 1\mskip-4.2mu l}
{\rm 1\mskip-4.6mu l}{\rm 1\mskip-5.2mu l}}
\newcommand\Ind[1]{{\One_{\{#1\}}}} 

\newcommand{\xra}{\xrightarrow} 

\newcommand{\set}[1]{\left\{#1\right\}}

\newcommand{\brc}[1]{\left(#1\right)}
\newcommand{\brcs}[1]{\left[#1\right]}

\renewcommand{\le}{\leqslant} 
\renewcommand{\ge}{\geqslant}

\usepackage{color}

\begin{document}

\begin{frontmatter}

\title{Minimax and pointwise sequential changepoint detection and identification for general stochastic models}


\author[A1]{Serguei M. Pergamenchtchikov}

\author[A2]{Alexander G. Tartakovsky\corref{mycorrespondingauthor}}

\author[A3]{Valentin S. Spivak}

\address[A1]{Laboratoire de Math\'ematiques Rapha\"el Salem, 
UMR 6085 CNRS-Universit\'e de Rouen Normandie, France and
National Research Tomsk State University, 
International Laboratory of Statistics of Stochastic Processes and Quantitative Finance, Tomsk, Russia}

\address[A2]{Space Informatics Laboratory, Moscow Institute of Physics and Technology, Dolgoprudny, Moscow Region, Russia
and AGT StatConsult, Los Angeles, California, USA}

\cortext[mycorrespondingauthor]{Corresponding author. Email address: \url{agt@phystech.edu}}

\address[A3]{Space Informatics Laboratory, Moscow Institute of Physics and Technology, Dolgoprudny, Moscow Region, Russia}

\begin{abstract}
This paper considers the problem of joint change detection and identification assuming multiple composite post-change hypotheses. We propose a multihypothesis changepoint detection-identification procedure that controls the probabilities of false alarm and wrong identification. We show that the proposed procedure is asymptotically minimax and pointwise optimal, minimizing moments of the detection delay as probabilities of false alarm and wrong identification approach zero. The asymptotic optimality properties hold for general stochastic models with dependent observations. We illustrate general results for detection-identification of changes in multistream Markov ergodic processes. We consider several examples, including an application to rapid detection-identification of COVID-19 in Italy. Our proposed sequential algorithm allows much faster detection of COVID-19 than standard methods. 
\end{abstract}


 
\end{frontmatter}



%
\section{Introduction}\label{sec:Intro}

As discussed in \cite{NikiforovIEEEIT95,TNB_book14,TartakovskyIEEEIT2021,Tartakovsky_book2020}, in a variety of applications it is important not only to quickly detect abrupt 
changes but also to diagnose them (e.g., to determine which change in a set of possible changes has occurred). This problem of change detection and diagnosis applies, for example, 
to rapid detection and identification of intrusions in computer networks, object detection with various sensors, integrity monitoring of navigation systems, and early detection and localization of epidemics. 
Often called Change Detection and Isolation, the problem is a generalization of the quickest change detection problem to the case of multiple post-change hypotheses and can be formulated as 
joint change detection and identification. Nikiforov~\cite{NikiforovIEEEIT95} first considered the change detection-isolation problem in a minimax setting for independent and identically distributed 
(i.i.d.) observations (in pre-change and post-change modes with different distributions) and simple 
post-change hypotheses. Several versions of the multihypothesis CUSUM-type and SR-type procedures, which have
minimax optimality properties in the classes of rules with constraints imposed on the average run length to a false alarm and conditional probabilities of false isolation, 
are proposed by Nikiforov~\cite{NikiforovIEEEIT00, NikiforovIEEEIT03}
and Tartakovsky~\cite{Tartakovsky-SQA08b}. Dayanik et al.~\cite{dayaniketal-AOR13} proposed an asymptotically optimal Bayesian detection-isolation rule assuming that the prior distribution 
of the change point is geometric also in the i.i.d. case. However, in many practical applications, the i.i.d.\ assumption is too restrictive -- the observations may be
either non-identically distributed or dependent or both, i.e., non-i.i.d. Also, the post-change distribution is usually not completely known. 
Lai~\cite{LaiIEEE00} provided a certain generalization for the non-i.i.d.\ case and composite hypotheses for a specific loss function. 
Recently, Tartakovsky~\cite{TartakovskyIEEEIT2021} developed a general asymptotic multistream Bayesian theory of sequential change detection and identification for low rates of false alarms and 
misidentification, assuming (1) there are multiple data streams, (2) the change occurs in some data stream(s) at an unknown random point in time, and (3) it is necessary to detect
the change as soon as possible and identify which data streams are affected. However, a non-Bayesian multistream change detection-identification theory for non-i.i.d.\ data is still missing.

The primary goal of this paper is to provide a general non-Bayesian asymptotic multistream change detection-identification theory (minimax and pointwise) for non-i.i.d.\ data and 
composite post-change hypotheses. 
This theory generalizes changepoint detection theory (with no identification) developed by Pergamenchtchikov and Tartakovsky~\cite{PerTar-JMVA2019}.  
In Section~\ref{sec:Model}, we describe the general stochastic model and provide basic notation. 
In Section~\ref{sec:Main-Cnds}, we introduce main conditions. In Section~\ref{sec:PF}, we introduce the change detection-identification rule. In Section~\ref{sec:LB}, we derive the information lower bounds 
for moments of the detection delay in the class of changepoint detection-identification rules with constraints imposed 
on the probabilities of false alarm and wrong identification.  
In Section~\ref{sec:Up-bnd}, we prove asymptotic optimality of the proposed detection-identification rule as the probabilities of false alarm and misidentification go to zero. We show that the lower bounds 
are attained for this procedure under very general conditions. In Section~\ref{sec:Mrk}, we illustrate general results for detection-identification of changes in Markov ergodic processes. 
In Section~\ref{sec:Ex}, we consider two examples -- detection-identification of changes in (1) the parameters of multivariate linear difference equations and (2) the correlation coefficients of multistream 
$p$-th order autoregressive models. In Section~\ref{sec:Epid}, we propose a specific model for epidemics and show that the proposed change detection-identification rule is asymptotically optimal. 
We also apply our rule for detection of COVID-19 in Italy and show that it allows for much earlier detection of COVID-19 than standard methods. 

\section{Basic notation} \label{sec:Model}

We consider the $N$ independent streams  of observations
$(X_\zs{1,l})_\zs{l\ge 1}\,,\ldots,(X_\zs{N,l})_\zs{l\ge 1}$.
For any $1\le i\le N$, $\nu\ge 0$ and $\theta_\zs{i}$ from an open set $\Theta_\zs{i}\subseteq\bbr^{m}$
we denote by  $\Pb_\zs{i,\nu,\theta_\zs{i}}$ the distribution of the observations $(X_\zs{i,l})_\zs{l\ge 1}$ in $\bbr^{\infty}$.
In the case when $\nu=\infty$, this distribution will be denoted by $\Pb^{*}$. 
We use the convention that $X_\zs{i,\nu}$ is the last pre-change observation. 
Write $\Xb^{n}_\zs{i}=(X_\zs{i,1},\dots,X_\zs{i,n})$ for the concatenation of the first $n$ observations in the $i$th data stream.
Let now for any $1\le i\le N$
\begin{equation}
\label{des-1}
 \left(f^{*}_\zs{i,j}(y_\zs{j}|y_\zs{1},\ldots,y_\zs{j-1})\right)_\zs{j\ge 1}
 \quad\mbox{and}\quad
 \left(f_\zs{i,\theta_\zs{i},j}(y_\zs{j}|y_\zs{1},\ldots,y_\zs{j-1})\right)_\zs{j\ge 1}
\end{equation} 
be  sequences of conditional densities of $X_\zs{i,j}$ 
given $\Xb^{j-1}_\zs{i}$  with respect to some non-degenerate $\sigma$-finite measure. 
Note that for $1\le i\le N$ the density $\q_\zs{i}$ of $\Xb^{n}_\zs{i}$ in $\bbr^{n}$ has the following form
\begin{equation}\label{dens-2}
\q_\zs{i,\nu,\theta_\zs{i}}(y_\zs{1},\ldots,y_\zs{n})=
\begin{cases}
 \q^{*}_\zs{i}(y_\zs{1},\ldots,y_\zs{n}) 
   & \quad \text{for}~~ \nu \ge n\,;
\\[4mm]
  \prod_{l=1}^{\nu}  f^{*}_\zs{i,l}(y_\zs{l}|y_\zs{1},\ldots,y_\zs{j-1}) \, \prod_{l=\nu+1}^{n}  
  f_\zs{i,\theta_\zs{i},l}(y_\zs{l}|y_\zs{1},\ldots,y_\zs{l-1})&  \quad \text{for}~~ \nu < n\,,
\end{cases}
\end{equation}
where $\q^{*}_\zs{i}(y_\zs{1},\ldots,y_\zs{n}) = \prod_\zs{l=1}^n f^{*}_\zs{i,l}(y_\zs{l}|y_\zs{1},\ldots,y_\zs{l-1})$.

Denote by $\gamma$ a random variable with values in $\{1,\ldots,N\}$
and assume that the change can occur only in  the data stream $(X_\zs{\gamma,l})_\zs{l\ge 1}$ with probability $\upsilon_\zs{i}=\P(\gamma=i)$. 
For $\nu\ge 1$ and 
$\theta=(\theta_\zs{1},\ldots,\theta_\zs{N})\in\Theta=\Theta_\zs{1}\times\ldots\times\Theta_\zs{N}$
 the joint density of the observations $\Xb^{n}_\zs{1},\ldots,\Xb^{n}_\zs{N}$ is given by
\begin{equation}
\label{total-den-1}
\p_\zs{\nu,\theta}(y_\zs{1,1},\dots,y_\zs{N,n})
=\sum^{N}_\zs{i=1}\upsilon_\zs{i}
\p_\zs{i,\nu,\theta}(y_\zs{1,1},\dots,y_\zs{N,n})
\,,
\end{equation}
where
\begin{equation}\label{jointforall}
\p_\zs{i,\nu,\theta_\zs{i}}(y_\zs{1,1},\dots,y_\zs{N,n})
=
\q_\zs{i,\nu,\theta_\zs{i}}(y_\zs{i,1},\ldots,y_\zs{i,n})\prod^{N}_\zs{l\neq i}\q^{*}_\zs{l}(y_\zs{l,1},\ldots,y_\zs{l,n})
\,.
\end{equation}

In the sequel we denote by $\cM$ the set of all Markov times with respect to the filtration $(\Fc_\zs{n})_\zs{n\ge 0}$ where
$\Fc_\zs{0}=\{\Omega,\varnothing\}$ and $\Fc_\zs{n}=\sigma\left\{X_\zs{i,j}\,,\,1\le i\le N\,,\,1\le j\le n \right\}$.

Note that when $n>k$ and $\vartheta\in\Theta_\zs{i}$ the Radon-Nykodim density (likelihood ratio)
\begin{equation}
\label{R-D-den-1}
\g^{*}_\zs{i,k,n}
=
\g^{*}_\zs{i,k,n}(\vartheta)
=
\frac{\d \Pb_\zs{i,k,\vartheta}}{\d \Pb^{*}}{\Big \vert}_\zs{\Fc_\zs{n}}
=e^{Z_\zs{i,n}^{k}(\vartheta)}\,,
\end{equation}
where
\begin{equation}\label{Znk_df}
Z_\zs{i,n}^{k} (\vartheta) = \sum_\zs{l=k+1}^{n}\, 
\log\,\frac{f_\zs{i,\vartheta,l}(X_\zs{i,l} \vert \Xb^{l-1}_\zs{i})}{f^{*}_\zs{i,l}(X_\zs{l}|\Xb^{l-1}_\zs{i})}
\end{equation}
is the log-likelihood ratio, and for any $(\theta_\zs{i},\theta_\zs{j})\in\Theta_\zs{i}\times \Theta_\zs{j}$
with $i\neq j$ the Radom-Nikodym density
\begin{equation}
\label{R-D-den-i->j}
\g_\zs{i,j,k,n}=
\g_\zs{i,j,k,n}(\theta_\zs{i}, \theta_\zs{j})
=
\frac{\d \Pb_\zs{i,k,\theta_\zs{i}}}{\d \Pb_\zs{j,k,\theta_\zs{j}}}{\Big \vert}_\zs{\Fc_\zs{n}}
=e^{Z_\zs{i,n}^{k}(\theta_\zs{i})-Z_\zs{j,n}^{k}(\theta_\zs{j})}\,.
\end{equation}

A sequential change detection-identification procedure $\delta$ is a pair $\delta=(T,d)$, where $T$ is a stoping time from $\cM$, i.e., 
for any $1\le i\le N$, $k\ge 0$ and $\vartheta\in\Theta_\zs{i}$
the probability $\Pb_\zs{i,k,\vartheta}(T<\infty)=1$, and $d$ is a decision rule, i.e., a random variable with the values
in $\{1,\ldots,N\}$ which is measurable with respect to the $\sigma$-field $\Fc_\zs{T}$.  
We denote by $\cS$ the class of all sequential procedures. For $r \ge 1$ and $\theta_i\in\Theta_\zs{i}$,
define the risk  for a sequential procedure $\delta=(T,d)\in\cS$ associated with the conditional $r$-th moment of the detection delay
\begin{equation}\label{SCrADD}
\Rc_\zs{i,k,\theta_i}(\delta)
=  
\EV_\zs{i,k,\theta_i}
\left[
 (T-k)^{r}\,\Chi_\zs{\{d=i\}}
 \vert T>k
 \right]
\,,
\end{equation}
where $\EV_\zs{i,k,\theta_i}$ is the expectation  with respect to the distribution 
$\Pb_\zs{i,k,\theta_i}$ in $\bbr^{\infty}$.

Introduce the conditional probability of false alarm $\Pb^{*}( T < k+\m^{*}\,,\,d=i \vert T\ge k)$ on the event $\{d=i\}$ in the interval  $[k, k+m^{*})$,
i.e., the probability of raising the alarm with the decision $d=i$ that there is a change in the $i$th stream when there is no change.  Also, introduce the misidentification probabilities 
$\Pb_\zs{i,k,\theta_i}(d=j \vert T>k)$, $i \neq j$, $i,j=1,\dots,N$.

 For any $N\times N$ matrix $\beta=(\beta_\zs{i,j})_\zs{1\le i,j\le N}$ with $0<\beta_\zs{i,j}<1$,  $\m^{*}\ge 1$ and $\k^{*}>\m^{*}$ we introduce the class of change detection-identification rules
\begin{align}
\nonumber
\Hc(\beta, \k^{*},\m^{*})&=\left\{\delta\in\cS: \sup_{1\le k\le \k^{*}-\m^{*}}\,\max_\zs{1\le i\le N} 
\frac{\Pb^{*}\left( T < k+\m^{*}\,,\,d=i \vert T\ge k \right)}{\beta_\zs{i,i}} \le 1\,, \right.
\\\label{sec:Prbf.4}
&\,
\left.\max_\zs{0\le k\le \k^{*}}\,\max_\zs{1\le i\le N}
 \sup_\zs{\theta\in\Theta_\zs{i}}
 \max_\zs{1\le j\neq i \le N}
 \frac{\Pb_\zs{i,k,\theta_i}\left(d=j \vert T>k \right)}{\beta_\zs{i,j}}
  \le 1
\right\}\,.
\end{align} 

Our goal is to find a sequential procedure asymptotically optimal in two problems in the class of detection-identification rules $\Hc(\beta, \k^{*},\m^{*})$: the pointwise minimization
\begin{equation}\label{sec:Prbf.5-0}
\inf_\zs{\delta\in \Hc(\beta, \k^{*},\m^{*})}\, \Rc_\zs{i,k,\theta_i}(\delta) \quad \text{for every}~ 
1\le i\le N\,, k\ge 0
~ \text{and} ~ \theta_i\in\Theta_\zs{i}
\end{equation}
and the minimax optimization 
\begin{equation}\label{sec:Prbf.5}
\inf_\zs{\delta\in \Hc(\beta, \k^{*},\m^{*})}\,\sup_\zs{k\ge 0}\,\Rc_\zs{i,k, \theta_i}(\delta)
\quad \text{for every}~ 
1\le i\le N
~ \text{and} ~ \theta_i\in\Theta_\zs{i}\,.
\end{equation}
\noindent The parameters $\k^{*}$ and $\m^{*}$  will be specified later.

\section{Main conditions}\label{sec:Main-Cnds}

For a fixed $\theta_i \in \Theta_\zs{i}$, we assume the following conditions for  the log-likelihood ratio (LLR) 
processes $(Z_\zs{i,n}^{k}(\theta_i))_\zs{n \ge k+1}$ introduced in \eqref{Znk_df} for $1\le i\le N$ and $\theta_i\in\Theta_\zs{i}$.

\noindent $(\A_\zs{1}$) {\em  For any $1\le i,j\le N$
there  are $\Theta_\zs{i}\times \Theta_\zs{j} \to \bbr_\zs{+}$ positive continuous functions $I_\zs{i,j}$ with 
\begin{equation}
\label{sec:Main-Cnds-Lb-1Ps-1}
0<\inf_\zs{(\theta_\zs{i},\theta_\zs{j})\in\Theta_\zs{i}\times \Theta_\zs{j}}I_\zs{i,j}(\theta_\zs{i},\theta_\zs{j})
\le \sup_\zs{(\theta_\zs{i},\theta_\zs{j})\in\Theta_\zs{i}\times \Theta_\zs{j}}I_\zs{i,j}(\theta_\zs{i},\theta_\zs{j})
<\infty
\end{equation}
 such that for any $k\ge 0$, $\varepsilon>0$
and $\theta=(\theta_\zs{1},\ldots \theta_\zs{N})\in \Theta_\zs{1}\times \ldots \times \Theta_\zs{N}$ 
\begin{equation}\label{sec:Cnd-2}
\lim_\zs{n\to\infty}\,
\max_\zs{1\le i\neq j\le N}\,
\Pb_\zs{i,k,\theta_\zs{i}}\left(Z_\zs{i,k+n}^{k}(\theta_\zs{i}) - Z_\zs{j,k+n}^{k}(\theta_\zs{j}) >(1+\varepsilon)I_\zs{i,j}(\theta_\zs{i}, \theta_\zs{j})n
 \right)=0 
\end{equation}
and
\begin{equation}\label{sec:Cnd-1}
\lim_\zs{n\to\infty}
\max_\zs{1\le i\le N}\,
\Pb_\zs{i,k,\theta_\zs{i}}\left(Z_\zs{i,k+n}^{k}(\theta_\zs{i})>(1+\varepsilon)\wt{I}_\zs{i}(\theta_\zs{i})n
 \right)=0\,,
\end{equation}
where $\wt{I}_\zs{i}(\theta_i)=I_\zs{i,i}(\theta_i,\theta_i)$ for $\theta_i\in\Theta_\zs{i}$.
}

In order to study asymptotic approximations to risks of the change detection-identification rule introduced below in Section~\ref{sec:PF} and for establishing its asymptotic optimality, 
we impose the following left-tail conditions:

\noindent $(\A_{2}(r)$) {\em
For any $1\le i,j\le N$
there  are $\Theta_\zs{i}\times \Theta_\zs{j} \to \bbr_\zs{+}$ positive continuous functions $I_\zs{i,j}$ with  the property \eqref{sec:Main-Cnds-Lb-1Ps-1}
such that for every $\theta=(\theta_\zs{1},\ldots,\theta_\zs{N})\in\Theta_\zs{1}\times\ldots\times\Theta_\zs{i}=\Theta$ and for any  $0<\vae<1$
\begin{equation}\label{rcompLeft-1}
\lim_\zs{\zeta\to 0}\,
\sup_\zs{\theta\in\Theta}
\max_\zs{1\le i \le N}\,
\sum_{n=1}^\infty \, n^{r-1} \, \sup_\zs{k \ge 0} \Pb_\zs{i,k,\theta_\zs{i}}\brc{
 \inf_\zs{\vert u-\theta_\zs{i}\vert<\zeta} Z_\zs{i,k+n}^{k}(u) < (1  - \varepsilon) \wt{I}_\zs{i}(\theta_\zs{i})n} <\infty
\end{equation}
and
\begin{equation}\label{rcompLeft-2}
\lim_\zs{\zeta\to 0}\,
\sup_\zs{\theta\in\Theta}
\max_\zs{1\le i\neq j \le N}\,
\sum_\zs{n=1}^\infty \, n^{r-1} \, \sup_\zs{0\le k\le \zeta n} \Pb_\zs{i,k,\theta_\zs{i}}\brc{
 \inf_\zs{\vert u-\theta_\zs{i}\vert<\zeta} Z_\zs{i,k+n}^{k}(u) - Z^{*}_\zs{j,k+n} < (1  - \varepsilon)\wt{I}_\zs{i}(\theta_\zs{i})n} <\infty\,,
\end{equation}
where $Z^{*}_\zs{j,n}=\max_\zs{0\le l\le n}\sup_\zs{z\in\Theta_\zs{j}}Z_\zs{j,k+n}^{l}(z)$. }

\begin{remark}
This is always true for i.i.d. data models with Kullback--Leibler informations given by $$
I_\zs{i,j}(\theta_\zs{i}, \theta_\zs{j})=
\int \log\left(\frac{f_\zs{i,\theta_\zs{i}}(x)}{f_\zs{j,\theta_\zs{j}}(x)}\right) f_\zs{i,\theta_\zs{i}}(x) \d\mu(x)
\quad\mbox{and}\quad
\wt{I}_\zs{i}(\theta_\zs{i})= 
\int \log \left(\frac{f_\zs{i,\theta_\zs{i}}(x)}{f^{*}_\zs{i}(x)}\right) f_\zs{i,\theta_\zs{i}}(x) \d\mu(x) 
$$
for $(\theta_\zs{i}, \theta_\zs{j})\in\Theta_\zs{i}\times \Theta_\zs{j}$.
\end{remark}

\section{Sequential change detection-identification procedure} \label{sec:PF}

First introduce weight distributions $(W_\zs{i}(\vartheta))_\zs{1\le i\le N}$, which  are probability measures on the sets $\Theta_\zs{i}$, i.e. $W_\zs{i}(\Theta_\zs{i})=1$ for any 
$1\le i\le N$. In what follows, we assume that $W_\zs{i}(\cdot)$  satisfy the following condition:

\noindent
$(\C_\zs{W})$
{\em For any $\varepsilon>0$ and any 
$\vartheta\in\Theta_\zs{i}$ the measure  $W_\zs{i}\{u\in\Theta_\zs{i}\,:\,\vert u-\vartheta\vert<\varepsilon \}>0$.}

Now,  for some fixed $0<\varrho<1$ we set
\begin{equation}\label{sec:Prbf.1}
\pi_\zs{k}
=\pi_\zs{k}(\varrho)
=\varrho\,\left(1-\varrho\right)^{k}\,, \quad k =0,1,2,\dots 
\end{equation}
\noindent
and using this distribution we set
\begin{equation}
\label{L-1}
\L_\zs{i,n}
=
\sum^{n-1}_\zs{k=0}\,\pi_\zs{k}\,
 \int_\zs{\Theta_\zs{i}}\, 
\g^{*}_\zs{i,k,n}(\vartheta)
  W_\zs{i}(\rm{d}\vartheta)
\quad\mbox{and}\quad
{\wh{\L}}_\zs{i,n}
=
\sum^{n-1}_\zs{k=0}\,\pi_\zs{k}\,
 \sup_\zs{\theta\in\Theta_\zs{i}}\, 
 \g^{*}_\zs{i,k,n}(\vartheta)
 \,.
\end{equation}
\noindent
Using these statistics we define the following random $N\times N$ matrix $\U_\zs{n}$ as
\begin{equation}
\label{ST-RdMatr}
<\U_\zs{n}>_\zs{i,j}
=\frac{\L_\zs{i,n}}{\wh{\L}_\zs{j,n}}
\quad \mbox{if $i\neq j$}  \quad \text{and} \quad
<\U_\zs{n}>_\zs{i,i}=\frac{\L_\zs{i,n}}{\sum_\zs{l\ge n}\pi_\zs{l}}
\,,
\end{equation}
where $<\U>_\zs{i,j}$ is the $(i,j)$th element of the matrix $\U$.
Finally, using this matrix we set
\begin{equation}\label{WSR-def}
T^{*}_\zs{i,A}=\inf\set{n \ge 1:
 \min_\zs{1\le j\le N}\,
 \frac{<\U_\zs{n}>_\zs{i,j}}{A_\zs{i,j}}
 \ge 1}
 \,,
\end{equation}
where  $A=(A_\zs{i,j})_\zs{1\le i,j\le N}$ is a $N\times N$ matrix with positive elements which will be specified later.  In the definitions of stopping times we  set $\inf\{\varnothing\}=+\infty$. The sequential 
change detection-identification procedure $\delta^{*}_\zs{A}=(T^{*}_\zs{A},d^{*}_\zs{A})$ that will be studied in this paper has the form
\begin{equation}
\label{seq-prs-1}
T^{*}_\zs{A}=\min_\zs{1\le i\le N}\,T^{*}_\zs{i,A}
\quad\mbox{and}\quad
d^{*}_\zs{A}=i
\quad\mbox{if}\quad
T^{*}_\zs{i,A}=T^{*}_\zs{A}\,.
\end{equation}
 If there are several numbers $i$ for which $T^{*}_\zs{i,A}=T^{*}_\zs{A}$ we can take arbitrary.
Note that, as we will see later in Proposition~\ref{Pr.sec:Up-bnd-1},  the condition $(\A^{*}_{2}(r)$) implies that $T^{*}_\zs{A}$ is a $\Pb_\zs{i,k,\theta_i}$-proper stopping time, that is, 
for any $1\le i\le N$, $k\ge 1$ and $\theta_i\in\Theta_\zs{i}$ 
\begin{equation}
\label{Prp-11-**}
\Pb_\zs{i,k,\theta_i}(T^{*}_\zs{A}<\infty)=1\,.
\end{equation}
Now, for any sequential procedure $\delta=(T,d)\in\cS$ we set
\begin{equation}\label{sec:Prbf.1-2-0}
\PFA_\zs{i}(\delta)= \sum_\zs{k=0}^\infty \pi_\zs{k} \Pb^{*}\left(T \le k\,,\,d=i\right)
\end{equation}
and
\begin{equation}\label{sec:Prbf.2-2-0}
\PMI_\zs{i,j}(\delta)= 
\sum_\zs{k=0}^\infty \pi_\zs{k} 
\sup_\zs{\theta_i\in\Theta_\zs{i}}\,
\Pb_\zs{i,k,\theta_i}\left(T>k\,,\,d= j\right)\,.
\end{equation}
\noindent 
For some $N\times N$ matrix $\alpha=(\alpha_\zs{i,j})$ with $0<\alpha_\zs{i,j}<1$ and
 some fixed $0<\varrho<1$, define the following {\em Bayesian class}:
\begin{align}\label{sec:Prbf.2}
\Delta(\alpha,\varrho)= \left\{\delta\in\cS:\,
\max\left(
 \max_\zs{1\le i\le N}\frac{\PFA_\zs{i}(\delta)}{\alpha_\zs{i,i}}
\,,\,
\max_\zs{ 1 \le i, j  \le N, i \neq j}
\frac{\PMI_\zs{i,j}(\delta))}{\alpha_\zs{i,j}}
\right)
 \le 1\right \}\
 \,.
\end{align}
\noindent
Next, for any arbitrary fixed matrix $\beta=(\beta_\zs{i,j})_\zs{1\le i,j\le N}$ and $0<\varrho<1$ introduce two matrices 
 $\alpha_\zs{1}=(\alpha^{(1)}_\zs{i,j})_\zs{1\le i,j\le N}$ and $\alpha_\zs{2}=(\alpha^{(2)}_\zs{i,j})_\zs{1\le i,j\le N}$ as
\begin{equation} 
\label{sec:alph-1-2-16-03-1}
\alpha^{(1)}_\zs{i,j}=\frac{\beta_\zs{i,i}(1-\varrho)^{\k^{*}}}{1+\tr\beta}\Chi_\zs{\{i= j\}}
+
\frac{\beta_\zs{i,j}\varrho(1-\varrho)^{\k^{*}}}{1+\tr\beta}\Chi_\zs{\{i\neq j\}}
\end{equation}
and
\begin{equation} 
\label{sec:alph-1-2-16-03-2}
\alpha^{(2)}_\zs{i,j}=
\left(\beta_\zs{i,i}+(1-\varrho)^{\m^{*}+1}\right)
\Chi_\zs{\{i= j\}}
+
\left(\beta_\zs{i,j}+(1-\varrho)^{\k^{*}+1}\right)\,
\Chi_\zs{\{i\neq j\}}
\,.
\end{equation}
\noindent

The following proposition compares classes \eqref{sec:Prbf.4} and \eqref{sec:Prbf.2}.

\begin{proposition}\label{Pr:Classes} 
For any matrix $\beta=(\beta_\zs{i,j})_\zs{1\le i,j\le N}$, $1\le \m^{*}<\k^{*}$ and $0<\varrho<1$ the following inclusions hold
\begin{equation}
\label{sec:Comp-16-03}
\Delta(\alpha_\zs{1},\varrho)
\subset 
\Hc(\beta, \k^{*},\m^{*})
\subset
\Delta(\alpha_\zs{2},\varrho)
\,. 
\end{equation}
\end{proposition}
\proof First note that if $\delta=(T,d)\in \Delta(\alpha_\zs{1},\varrho)$, then for any $1\le i\le N$  and $k\ge 1$
$$
\alpha^{(1)}_\zs{i,i}\ge 
\sum_\zs{l=0}^\infty \pi_\zs{l} \Pb^{*}\left(T \le l\,,\,d=i\right)
\ge \Pb^{*}\left(T \le k\,,\,d=i\right) \sum_\zs{l=k}^\infty \pi_\zs{l} 
=(1-\varrho)^{k}\Pb^{*}\left(T \le k\,,\,d=i\right)\,,
$$
i.e., for $k\ge 1$
\begin{equation} 
\label{Up-tau-1}
\Pb^{*}\left(T \le k\,,\,d=i\right)\le (1-\varrho)^{-k}\alpha^{(1)}_\zs{i,i}
\quad\mbox{and}\quad
\Pb^{*}\left(T \le k\right)\le 
\sum^{N}_\zs{j=1}
\Pb^{*}\left(T \le k\,,\,d=j\right)\le (1-\varrho)^{-k}
\tr\alpha_\zs{1}
\,.
\end{equation}
\noindent 
Therefore, for any $1\le i\le N$ and $1\le k\le \k^{*}-\m^{*}$
\begin{align*}
\Pb^{*}&\left(T < k+\m^{*}\,,\,d=i \vert  T\ge k \right)
=\frac{\Pb^{*}\left(k\le T < k+\m^{*}\,,\,d=i\right)}{1-\Pb^{*}(T < k)}\\[2mm]
&\le 
\frac{\Pb^{*}\left( T \le \k^{*}\,,\,d=i\right)}{1- \Pb^{*}(T \le \k^{*})}
\le 
\frac{(1-\varrho)^{-\k^{*}}\alpha^{(1)}_\zs{i,i}}{1- (1-\varrho)^{-\k^{*}} \tr \alpha_\zs{1}}
:=\beta_\zs{i,i}
\,.
\end{align*}
\noindent
Moreover, for $i\neq j$ and any $k\ge 1$
$$
\alpha^{(1)}_\zs{i,j} \ge 
\pi_\zs{k} 
\sup_\zs{\theta_i\in\Theta_\zs{i}}\,
\Pb_\zs{i,k,\theta_i}\left(T>k\,,\,d= j\right)
=\varrho(1-\varrho)^{k}
\sup_\zs{\theta_i\in\Theta_\zs{i}}\,
\Pb_\zs{i,k,\theta_i}\left(T>k\,,\,d= j\right)\,,
$$
i.e., in view of \eqref{Up-tau-1}, for $1\le k\le \k^{*}$ and $\theta_i\in\Theta_\zs{i}$
\begin{align*}
\Pb_\zs{i,k,\theta_i}&\left(d= j \vert T>k\right)
=
\frac{\Pb_\zs{i,k,\theta_i}\left(T >k\,,\,d= j\right)}{\Pb_\zs{i,k,\theta_i}\left(T>k\right)}
=
\frac{\Pb_\zs{i,k,\theta_i}\left(T>k\,,\,d= j\right)}{1-\Pb^{*}\left(T>k\right)}
\\[2mm]
&\le 
\frac{
\alpha^{(1)}_\zs{i,j}\varrho^{-1}(1-\varrho)^{-k}}{1-\Pb^{*}\left(T\le k\right)}
\le 
\frac{
\alpha^{(1)}_\zs{i,j}\varrho^{-1}(1-\varrho)^{-\k^{*}}}{1- (1-\varrho)^{-\k^{*}} \tr \alpha_\zs{1}}
=\beta_\zs{i,j}
\,.
\end{align*}
This implies, that $\delta\in \Hc(\beta, \k^{*},\m^{*})$, i.e., we get the first inclusion in  \eqref{sec:Comp-16-03}. 

Let now
$\delta=(T,d)\in \Hc(\beta, \k^{*},\m^{*})$, i.e., for any $1\le k\le \k^{*}-\m^{*}$ and $1\le i\le N$
$$
\beta_\zs{i,i}\ge 
\Pb^{*}\left( T < k+\m^{*}\,,\,d=i\vert T\ge k \right) 
=\frac{\Pb^{*}\left(k\le T < k+\m^{*}\,,\,d=i\right) }{\Pb^{*}(T\ge k) }
\ge \Pb^{*}\left(k\le T < k+\m^{*}\,,\,d=i\right) 
$$
and, in particular,
$$
 \beta_\zs{i,i}\ge
\Pb^{*}\left(T < 1+\m^{*}\,,\,d=i \right)\,.
$$
Therefore,
\begin{align*}
\sum_\zs{k\ge 0}\pi_\zs{k}\Pb^{*}\left(T \le k\,,\,d=i \right)
&\le \Pb^{*}\left(T < 1+\m^{*}\,,\,d=i \right) + \sum_\zs{k\ge \m^{*}+1}\pi_\zs{k}\\[2mm]
&\le \beta_\zs{i,i}+(1-\varrho)^{\m^{*}+1} =\alpha^{(2)}_\zs{i,i}\,.
\end{align*}
Furthermore, for any $i\neq j$, $\theta_i\in\Theta_\zs{i}$ and $0\le k\le \k^{*} $
$$
\beta_\zs{i,j}\ge
\Pb_\zs{i,k,\theta_i}\left(d=j \vert T>k \right)
 \ge
 \Pb_\zs{i,k,\theta_i}\left(k<T <\infty\,,\,d=j \right), 
$$
i.e.,
$$
\sum_\zs{k\ge 0}\pi_\zs{k} \Pb_\zs{i,k,\theta_i}\left(T>k\,,\,d=j \right)
\le \beta_\zs{i,j} + \sum_\zs{k> \k^{*}}\pi_\zs{k} =\beta_\zs{i,j}+(1-\varrho)^{\k^{*}+1}=\alpha_\zs{i,j}\,.
$$
Thus, we obtain the last inclusion \eqref{sec:Comp-16-03}.
 Hence Proposition \ref{Pr:Classes}.
\endproof

The first question we ask is how to select the thresholds in the procedure \eqref{seq-prs-1} to imbed it into class $\Delta(\alpha,\rho)$. 
To study this question we need the following probability measures on $\bbn\times \bbr^{\infty}$ which for any $1\le i\le N$ are defined  as
\begin{equation}
\label{sec:Baes-measur-1}
\Q_\zs{i}(J\times A)=\sum_\zs{k\in J}
\pi_\zs{k}\,
\int_\zs{\Theta_\zs{i}}
\Pb_\zs{i,k,\vartheta}(A)\,W_\zs{i}(\d \vartheta)
\,,\quad J\subseteq \bbn\quad\mbox{and}\quad A\in\cB(\bbr^{\infty})\,,
\end{equation}
where $\bbn=\{0,1,\ldots\}$ and $\cB(\bbr^{\infty})$ is the cylinder field in $\bbr^{\infty}$.  In the sequel we denote by
$\E^{\Q_\zs{i}}$ the expectation over the probability measure $\Q_\zs{i}$. 
One can check directly that
\begin{equation}
\label{sec:Cond-Expct-1}
\Q_\zs{i}\left(\nu\ge n\vert \cF_\zs{n} \right)
=\frac{1}{1+<\U_\zs{n}>_\zs{i,i}}
\,,
\end{equation}
where the matrix $\U_\zs{n}$ is defined in \eqref{ST-RdMatr}.
 
\begin{lemma}\label{Lem:PFASR} 
For all $1\le i\le N$  the probabilities  \eqref{sec:Prbf.1-2-0} satisfy the inequalities 
\begin{equation}
\label{upper-bnd-1}
\PFA_\zs{i}(\delta^{*}_\zs{A}) \le \frac{1}{1+ A_\zs{i,i}}\,. 
\end{equation}
\end{lemma}

\proof
Note that
\begin{align*}
\PFA_\zs{i}(\delta^{*}_\zs{A})&= \sum_\zs{k=0}^\infty \pi_\zs{k} \Pb^{*}\left(T^{*}_\zs{A} \le k\,,\,d=i\right)=
 \sum_\zs{k=0}^\infty \pi_\zs{k}\int_\zs{\Theta_\zs{i}} \Pb_\zs{i,k, \vartheta}\left(T^{*}_\zs{i,A} \le k\,,\,d=i\right)\,W_\zs{i}(\d \vartheta)\\[2mm]
&\le
\sum_\zs{k=0}^\infty \pi_\zs{k}\int_\zs{\Theta_\zs{i}} \Pb_\zs{i,k,\vartheta}\left(T^{*}_\zs{i,A} \le k\right)\,W_\zs{i}(\d \vartheta)
=\Q_\zs{i}\left(\nu \ge T^{*}_\zs{i,A} \right)
= 
\E^{\Q_\zs{i}}
\left[\Q_\zs{i}\left(\nu \ge T^{*}_\zs{i,A} \vert \cF_\zs{T^{*}_\zs{i,A}} \right)\right]\,.
\end{align*}
\noindent 
Therefore, using \eqref{sec:Cond-Expct-1} and
\eqref{WSR-def}
we obtain that
$$
\PFA_\zs{i}(\delta^{*}_\zs{A})\le 
\E^{\Q_\zs{i}}\brcs{\frac{1}{1+<\U_\zs{T^{*}_\zs{i,A}}>_\zs{i,i}}}
\le \frac{1}{1+A_\zs{i,i}}\, ,
$$
which completes the proof.
 \endproof

\begin{lemma}\label{Lem:PFASR-1} 
For any $1\le i, j\le N$, $i \neq j$  the PMI probabilities of the procedure  \eqref{seq-prs-1}
satisfy the inequalities
\begin{equation}
\label{upper-bnd-2}
\PMI_\zs{i,j}(\delta^{*}_\zs{A})
\le \frac{1}{A_\zs{j,i}}\,.
\end{equation}
\end{lemma}

\proof
\noindent
First note that  for the  rule \eqref{WSR-def} we obtain that  for any $j\neq i$ and $\theta_i\in\Theta_\zs{i}$ 
\begin{align*}
\Pb_\zs{i,k,\theta_i}\left(T^{*}_\zs{A}>k\,,\,d=j\right)&\le 
\Pb_\zs{i,k,\theta_i}\left(k<T^{*}_\zs{j,A} <\infty\right)\\
&\le 
\frac{1}{A_\zs{j,i}}\,\EV_\zs{i,k,\theta_i} \brcs{<\U_\zs{T^{*}_\zs{j,A}}>_\zs{j,i} \Chi_\zs{\left\{k<T^{*}_\zs{j,A} <\infty\right\}}}\\
&=
\frac{1}{A_\zs{j,i}}\,\EV^{*} \brcs{<\U_\zs{T^{*}_\zs{j,A}}>_\zs{j,i}\g_\zs{i,k,T^{*}_\zs{j,A}}(\theta_i) \Chi_\zs{\left\{k<T^{*}_\zs{j,A} <\infty\right\}}}
\,.
\end{align*}
\noindent
Therefore, in view of the definition in \eqref{L-1}, we get
\begin{align*}
\sum_\zs{k\ge 0}\pi_\zs{k}
\Pb_\zs{i,k,\theta_i}\left(T^{*}_\zs{A}>k\,,\,d=j\right)
&\le 
\frac{1}{A_\zs{j,i}}\,\EV^{*} \brcs{<\U_\zs{T^{*}_\zs{j,A}}>_\zs{j,i}
\sum^{T^{*}_\zs{j,A}-1}_\zs{k=0}\pi_\zs{k}
\g_\zs{i,k,T^{*}_\zs{j,A}}(\theta_i)
 \Chi_\zs{\left\{T^{*}_\zs{j,A} <\infty\right\}}}\\
&\le 
\frac{1}{A_\zs{j,i}}\,\EV^{*} \brcs{<\U_\zs{T^{*}_\zs{j,A}}>_\zs{j,i}
\wh{\L}_\zs{i,T^{*}_\zs{j,A}}\,
 \Chi_\zs{\left\{T^{*}_\zs{j,A} <\infty\right\}}}\\
 &=
 \frac{1}{A_\zs{j,i}}\,\EV^{*} \brcs{
\L_\zs{j,T^{*}_\zs{j,A}}\,
 \Chi_\zs{\left\{T^{*}_\zs{j,A} <\infty\right\}}}
 \,.
\end{align*}
Moreover, note that
\begin{align*}
\EV^{*} \brcs{
\L_\zs{j,T^{*}_\zs{j,A}}\,
 \Chi_\zs{\left\{T^{*}_\zs{j,A} <\infty\right\}}}
&=\sum_\zs{k\ge 0}\pi_\zs{k}
 \int_\zs{\Theta_\zs{i}}\, 
\EV^{*}\brcs{
\g_\zs{j,k,T^{*}_\zs{j,A}}(\theta_i)\Chi_\zs{\left\{k<T^{*}_\zs{j,A} <\infty\right\}}}
  W_\zs{i}(\rm{d}\theta_i)\\[2mm]
&=
\sum_\zs{k\ge 0}\pi_\zs{k}
 \int_\zs{\Theta_\zs{i}}\, 
 \Pb_\zs{i,k,\theta_i}
\left(k<T^{*}_\zs{j,A} <\infty\right)\,
  W_\zs{i}(\rm{d}\theta_i))
  \le 1\,,
\end{align*}
which implies upper bound \eqref{upper-bnd-2}
\endproof

Now, if we take in \eqref{WSR-def}
\begin{equation}
\label{upper-bnd-2new}
A_\zs{i,j}
=
\left(\frac{1}{\alpha_\zs{i,i}}-1\right)
\Chi_\zs{\{i= j\}}
+
\frac{1}{\alpha_\zs{j,i}}\Chi_\zs{\{i\neq j\}}
\end{equation}
then using the property
\eqref{Prp-11-**} and the upper bounds \eqref{upper-bnd-1} and \eqref{upper-bnd-2}
we obtain that under condition $(\A^{*}_{2}(r)$)  the sequential procedure \eqref{seq-prs-1} belongs to class $\Delta(\varrho,\alpha)$ for any $0<\varrho<1$
and $\alpha=(\alpha_\zs{i,j})_\zs{1\le i,j\le N}$ with $0<\alpha_\zs{i,j}<1$.
Therefore, if we take 
\begin{equation}
\label{upper-bnd-2-1**}
A_\zs{i,j}
=
\left(\frac{1+\tr \beta}{\beta_\zs{i,i}(1-\varrho)^{\k^{*}}}-1\right)
\Chi_\zs{\{i= j\}}
+
\frac{1+\tr \beta}{\beta_\zs{j,i}\varrho (1-\varrho)^{\k^{*}}}
\Chi_\zs{\{i\neq j\}}
\,,
\end{equation}
we obtain that for any $0<\varrho<1$  the sequential procedure \eqref{seq-prs-1} belongs to class $\Hc(\beta, \k^{*},\m^{*})$.

\section{Information lower bounds}\label{sec:LB}

\subsection{Bayesian setting}\label{sec:Bay} 

For any matrix $\alpha=(\alpha_\zs{i,j})_\zs{1\le i,j\le N}$ and
any  parameter value $\theta_i\in\Theta_\zs{i}$ define
\begin{equation}
\label{Lower-bound-1}
\b_\zs{i,\alpha}(\theta_i)=\max_\zs{1\le j\le N}\,
\frac{\vert\log \alpha_\zs{j,i}\vert}{\iota_\zs{i,j}(\theta_i)}
\quad\mbox{and}\quad
\iota_\zs{i,j}(\theta_i)=\wt{I}_\zs{i}(\theta_i)\,\Chi_\zs{\{j=i\}}
+\wh{I}_\zs{i,j}(\theta_i)\,\Chi_\zs{\{j\neq i\}}\,,
\end{equation}
where the function $\wt{I}_\zs{i}(\cdot)$ is defined in
\eqref{sec:Cnd-1} and   
$\wh{I}_\zs{i,j}(\theta_i)=\inf_\zs{\vartheta\in\Theta_\zs{j}} I_\zs{i,j}(\theta_i, \vartheta)$ for $\theta_i\in\Theta_\zs{i}$.
In what follows we will always suppose without special emphasis that
\[
\inf_\zs{\theta_j\in\Theta_\zs{j}} I_\zs{i,j}(\theta_i, \theta_j)>0 \quad \text{for all $\theta_i\in\Theta_i$, $j \neq i$ and $i=1,\dots,N$}
\]
(see condition  \eqref{sec:Main-Cnds-Lb-1Ps-1}).

Write $\alpha_\zs{max}=\max_\zs{1\le i,j\le N}\alpha_\zs{i,j}$ and $\Delta_\alpha=\Delta(\alpha, \varrho_\alpha)$ (in case where $\varrho=\varrho_\alpha$ depends on $\alpha$). 
The following theorem establishes information lower bounds  in the Bayesian problem. These bonds will be used to obtain asymptotic lower bounds for 
$\Rc_{i,k,\theta_i}(\delta)$ in class $\Hc\left(\beta,\k^{*},\m^{*}\right)$ (see Theorem~\ref{Th.sec:Cnrsk.1}) and to prove asymptotic optimality of the proposed 
detection-identification procedure in this class.

\begin{theorem} \label{Th.sec:Bay.1} 
 Assume that the right-tail probability convergence condition $(\A_\zs{1})$  holds and in \eqref{sec:Prbf.1} the parameter of the geometric prior distribution
 $\varrho$ is a function of $\alpha$, i.e., $\varrho=\varrho_\zs{\alpha}$, such that
\begin{equation}
\label{cnd-alpha-1}
\lim_\zs{\alpha_\zs{max}\to 0}\left(\varrho_\zs{\alpha}+\frac{\vert\log \varrho_\zs{\alpha}\vert}{\vert \log \alpha_\zs{max}\vert}
\right)=0
\,. 
\end{equation}
Then, for any $r\ge 1$, $k \ge 0$, $\theta_i\in\Theta_\zs{i}$, and $1\le i\le N$ the following asymptotic lower bounds hold: 
\begin{equation} \label{sec:Bay.1-nw-1} 
\liminf_\zs{\alpha_\zs{max}\to 0}  \frac{\inf_\zs{\delta\in \Delta_\zs{\alpha}} 
\,
\EV_\zs{i,k,\theta_i} \brcs{(T-k)^{r}_\zs{+}\Chi_\zs{\{d=i\}}}
}{|\b_\zs{i,\alpha}(\theta_i)|^{r}}\, 
 \ge 1\,.
\end{equation}
\end{theorem}

\proof
To prove this theorem it suffices to show that for any $j\neq i$ and 
$(\theta_\zs{i},\theta_\zs{j})\in\Theta_\zs{i}\times \Theta_\zs{j}$ ($j=1,\dots,N$)
\begin{equation}
\label{LW-bnd-jneqi}
\inf_\zs{\delta\in \Delta_\zs{\alpha}}\,
\EV_\zs{i,k,\theta_\zs{i}}\,\brcs{(T-k)^{r}_\zs{+}\Chi_\zs{\{d=i\}}}
\ge
(1+\psi_\zs{\alpha,i,j}(\theta_\zs{i},\theta_\zs{j}))
\frac{\vert\log \alpha_\zs{j,i}\vert^{r}}{I^{r}_\zs{i,j}(\theta_i,\theta_j)}
\,,
\end{equation}
and for any $j=i$ and $\theta_i\in\Theta_i$
\begin{equation}
\label{LW-bnd-jeqi}
\inf_\zs{\delta\in \Delta_\zs{\alpha}}\,
\EV_\zs{i,k,\theta_\zs{i}}\,\brcs{(T-k)^{r}_\zs{+}\Chi_\zs{\{d=i\}}}
\ge
(1+\psi_\zs{\alpha,i,i}(\theta_\zs{i},\theta_\zs{i}))
\frac{\vert\log \alpha_\zs{i,i}\vert^{r}}{\widetilde{I}^{r}_\zs{i}(\theta_i)}
\,,
\end{equation}
where the term $\psi_\zs{\alpha,i,j}(\theta_\zs{i},\theta_\zs{j})$ is such that 
$$
\lim_\zs{\alpha_\zs{max}\to 0}\,\vert \psi_\zs{\alpha,i,j}(\theta_\zs{i},\theta_\zs{j})\vert =0 \quad \text{for any 
$(\theta_\zs{i},\theta_\zs{j})\in\Theta_\zs{i}\times \Theta_\zs{j}$ and $1\le j\le N$}
\,.
$$
To prove \eqref{LW-bnd-jneqi} note that condition \eqref{sec:Cnd-2} implies  that
for any $\vae>0$, $k\ge 0$ and $(\theta_\zs{i},\theta_\zs{j}) \in\Theta_\zs{i}\times  \Theta_\zs{j}$ with $j\neq i$
\begin{equation} \label{sec:Bay.2-0}
\Pb_\zs{i,k,\theta_i}\set{Z_\zs{i,k+M}^{k}(\theta_\zs{i})-Z_\zs{j,k+M}^{k}(\theta_\zs{j})\ge (1+\vae) I_\zs{i,j} M}\xra[M\to\infty]{}0\,,
\end{equation}
where $I_\zs{i,j}=I_\zs{i,j}(\theta_\zs{i}, \theta_\zs{j})$.
\noindent
Define for $(\theta_\zs{i},\theta_\zs{j}) \in\Theta_\zs{i}\times  \Theta_\zs{j}$
$$
\D_\zs{i,j,k}(\delta)=\Pb_\zs{i,k,\theta_\zs{i}}
\left(k< T\le k+M_\zs{i,j}\,,\,d=i\right)
\quad\mbox{and}\quad
M_\zs{i,j}=M_\zs{i,j}(\theta_\zs{i},\theta_\zs{j})= (1-\varepsilon)\frac{|\log\alpha_\zs{j,i}|}{I_\zs{i,j}}
\,.
$$
 We now show that for any $k \ge 0$, $0<\varepsilon<1$ and $1\le j \neq i \le N$
\begin{equation} \label{sec:Bay.3}
\lim_\zs{\alpha_\zs{max}\to 0}\sup_\zs{\delta\in \Delta_\zs{\alpha}}\,
\D_\zs{i,j,k}(\delta)
\,=0\,.
\end{equation}
Using definition \eqref{R-D-den-i->j} we can  obtain that for $\m=k+M_\zs{i,j}$
\begin{equation} \label{sec:Bay.4}
\begin{aligned} 
 \D_\zs{i,j,k}(\delta)
 &=
 \EV_\zs{j,k,\theta_\zs{j}}\,\brcs{ \g_\zs{i,j,k,\m}\,\Chi_\zs{\{k\le T\le \m\,,\,d=i\}}}
  \le 
e^{(1+\varepsilon)I_\zs{i,j}  M_\zs{i,j}}\,
\Pb_\zs{j,k,\theta_\zs{j}}\left(k< T\le \m\,,\,d=i\right)\\[3mm]
 &+
\Pb_\zs{i,k,\theta_\zs{i}} 
 \left(Z_\zs{i,\m}^{k}(\theta_\zs{i})-Z_\zs{j,\m}^{k}(\theta_\zs{j})\ge (1+\varepsilon)I_\zs{i,j} M_\zs{i,j} \right)\,.
\end{aligned}
\end{equation}
Using the definition of $\PMI_{ji}(\delta)$ in \eqref{sec:Prbf.2-2-0} along with the fact that $\PMI_{ji}(\delta) \le \alpha_{j,i}$ for any $\delta\in \Delta(\alpha,\varrho_\alpha)$ we get 
\begin{align*}
\alpha_{j,i} &\ge \sum_{l \ge 0}{\pi_l} \sup_{\theta_j \in \Theta_j} \Pb_{j,l,\theta_j} (T>l, d=i) \ge \pi_k \Pb_{j,l,\theta_j} (T>l, d=i) \ge \pi_k  \Pb_{j,k,\theta_j} (k < T \le \m, d=i) 
\\
& = \varrho_\alpha (1-\varrho_\alpha)^k \Pb_{j,k,\theta_j} (k < T \le \m, d=i) \quad \text{for any $\theta_j\in\Theta_j$ and any $k \ge 0$},
\end{align*}
so that
$$
\sup_\zs{\delta\in\Delta_\zs{\alpha}}\, \Pb_\zs{j,k,\theta_\zs{j}}\left(k< T\le \m\,,\,d=i\right)
 \le \varrho_\zs{\alpha}^{-1}(1-\varrho_\zs{\alpha})^{-k} \alpha_\zs{j,i}
 =e^{-\vert\log \alpha_\zs{j,i}\vert-\log \varrho_\zs{\alpha} + k\varpi_\zs{\alpha}} 
 \,,
$$
where in view of \eqref{cnd-alpha-1} the term $\varpi_\zs{\alpha}=-\log(1-\varrho_\zs{\alpha})\to 0$ as $\alpha_\zs{max}\to 0$.
So the first term on the right-hand side of the inequality  \eqref{sec:Bay.4} can be estimated as
$$
\exp\left\{(1+\varepsilon)I_\zs{i,j}  M_\zs{i,j}
-\log \varrho_\zs{\alpha} + k\varpi_\zs{\alpha}+\log \alpha_\zs{j,i}
\right\}
\le 
\exp\left\{-\varepsilon^{2}\vert\log \alpha_\zs{j,i}\vert
-\log \varrho_\zs{\alpha} + k\varpi_\zs{\alpha}
\right\}
\le 
\,\exp\set{-\varepsilon^{2}\vert\log \alpha_\zs{max}\vert
-\log \varrho_\zs{\alpha} + k\varpi_\zs{\alpha}}
$$
and by condition \eqref{cnd-alpha-1} it goes to zero as $\alpha_\zs{max}\to 0$.  Therefore, \eqref{sec:Bay.4} and \eqref{sec:Bay.2-0}
impliy \eqref{sec:Bay.3} for any $j\neq i$.

 Let now $i=j$. Using the definition \eqref{R-D-den-1} we can rewrite the inequality \eqref{sec:Bay.4} as
\begin{equation}
\label{Up-Bnd-11-2}
\begin{aligned} 
 \D_\zs{i,i,k}(\delta)
 &=
 \EV^{*}\, \brcs{\g^{*}_\zs{i,k,\m}\,\Chi_\zs{\{k\le T\le \m\,,\,d=i\}}}
  \le 
e^{(1+\varepsilon)\wt{I}_\zs{i}(\theta_i) M_\zs{i,i}}\,
\Pb^{*}
\left(k< T\le \m\,,\,d=i\right)\\
 &+
\Pb_\zs{i,k,\theta_\zs{i}} \left(Z_\zs{i,\m}^{k}(\theta_\zs{i})
 \ge (1+\varepsilon)\wt{I}_\zs{i}(\theta_i) M_\zs{i,i} \right)\,,
\end{aligned}
\end{equation}
where $\wt{I}_\zs{i}(\theta_i)=I_\zs{i,i}(\theta_i,\theta_i)$ and where by condition \eqref{sec:Cnd-1}   
\begin{equation}\label{sec:Cnd-1-11}
\lim_\zs{M_{ii}\to\infty}
\Pb_\zs{i,k,\theta_\zs{i}}\left(Z_\zs{i,k+M_{ii}}^{k}(\theta_\zs{i})>(1+\varepsilon)\wt{I}_\zs{i}(\theta_\zs{i})M_{ii}
 \right)=0 \quad \text{for any $\theta_i\in\Theta_\zs{i}$}\,.
\end{equation}
Now, the definition of class $\Delta_\alpha=\Delta(\alpha,\varrho_\alpha)$ in \eqref{sec:Prbf.2} implies that for any $\delta \in \Delta_\alpha$, any $k\ge 1$ and all $i =1,\dots,N$
\begin{align*}
\alpha_{ii} \ge \sum_{\ell \ge 1} \pi_\ell \Pb^*(T \le l, d=i) \ge \Pb^*(T \le k , d=i) \sum_{\ell =1}^k \varrho_\alpha(1-\varrho_\alpha)^{l-1} = (1-\varrho_\alpha)^{k} \, \Pb^*(T \le k , d=i) ,
\end{align*}
which yields
\begin{equation}\label{sec:Bay.5}
\sup_\zs{\delta\in\Delta_\zs{\alpha}}\,
\Pb^{*}\left(T\le k\,,\,d=i\right)\le \alpha_\zs{i,i}\,(1-\varrho_\zs{\alpha})^{-k}=e^{-\vert\log\alpha_\zs{i,i}\vert+k\varpi_\zs{\alpha}}\,. 
\end{equation}

Therefore, the first term on the right side of the inequality \eqref{Up-Bnd-11-2}
may be estimated as
$$
e^{(1+\varepsilon)\wt{I}_\zs{i,i}(\theta_i)  M_\zs{i,i}-|\log \alpha_\zs{i,i}|
+\varpi_\zs{\alpha}k+\varpi_\zs{\alpha}M_\zs{i,i}}
\le e^{-\varepsilon^{2} |\log\alpha_\zs{i,i}|+\varpi_\zs{\alpha}k+\varpi_\zs{\alpha}M_\zs{i,i}} 
$$
and it goes to zero for any fixed $0\le k<\infty$ as $\alpha_\zs{max}\to 0$, which along with  \eqref{sec:Cnd-1-11} implies \eqref{sec:Bay.3} for $i=j$. 
 
To obtain lower bounds \eqref{LW-bnd-jneqi} and \eqref{LW-bnd-jeqi} note that for any  $1\le j\le N$ 
\begin{align}\nonumber
\EV_\zs{i,k,\theta_\zs{i}}\,\brcs{(T-k)^{r}_\zs{+}\Chi_\zs{\{d=i\}}}
&\ge
\EV_\zs{i,k,\theta_\zs{i}}\,\brcs{(T-k)^{r}_\zs{+}\Chi_\zs{\{T > k+M_\zs{i,j}\,,\,d=i\}}}
\ge M^{r}_\zs{i,j}\,
\Pb_\zs{i,k,\theta_\zs{i}}\left( T > k+M_\zs{i,j}\,,\,d=i\right)
\\ \label{sec:01-05-1}
&
=
M^{r}_\zs{i,j}\left(
\Pb_\zs{i,k,\theta_\zs{i}}
\left(T > k\,,\,d=i\right) -\D_\zs{i,j,k}(\delta)\right)\,.
\end{align}
\noindent Using the upper bound  \eqref{sec:Bay.5}, we get
\begin{align*}
\Pb_\zs{i,k,\theta_\zs{i}}
\left(T > k\,,\,d=i\right)&=
\Pb_\zs{i,k,\theta_\zs{i}}
\left(d=i\right)
-\Pb_\zs{i,k,\theta_\zs{i}}
\left(T \le k\,,\,d=i\right)\\
&=
\Pb_\zs{i,k,\theta_\zs{i}}
\left(d=i\right)
-
\Pb^{*}
\left(T \le k\,,\,d=i\right)\\
&\ge 1-\sum^{N}_\zs{l=1,l\neq i}
\Pb_\zs{i,k,\theta_\zs{i}}
\left(d= l\right)
-\alpha_\zs{i,i} (1-\varrho_\zs{\alpha})^{-k}
\,.
\end{align*}
Next, it follows from 
\eqref{sec:Prbf.2-2-0} and \eqref{sec:Prbf.2}
that for $l\neq i$
$$
\sup_\zs{\theta_i\in\Theta_\zs{i}}
\Pb_\zs{i,k,\theta_i}\left(T > k\,,\,d= l\right)
\le \frac{\alpha_\zs{i,l}}{\varrho_\zs{\alpha}}(1-\varrho_\zs{\alpha})^{-k}
=e^{\log \alpha_\zs{i,l}-\log \varrho_\zs{\alpha}+k\varpi_\zs{\alpha}}
\,,
$$
i.e., for any $l\neq i$ and $\theta_i\in\Theta_\zs{i}$
\begin{align*}
\Pb_\zs{i,k,\theta_i}
\left(d= l\right) &=
\Pb^{*}\left(T\le k\,,\,d= l\right)
+
\Pb_\zs{i,k,\theta_i}\left(T>k\,,\,d=l\right)\\[2mm]
&\le 
\alpha_\zs{l,l} (1-\varrho)^{-k}
+
\frac{\alpha_\zs{i,l}}{\varrho_\zs{\alpha}}(1-\varrho_\zs{\alpha})^{-k}
=
e^{\log \alpha_\zs{i,l}+k\varpi_\zs{\alpha}}
+
e^{\log \alpha_\zs{i,l}-\log \varrho_\zs{\alpha}+k\varpi_\zs{\alpha}} .
\end{align*}
Thus, in view of \eqref{cnd-alpha-1} for any $\theta_i\in\Theta_\zs{i}$
$$
\Pb_\zs{i,k,\theta_i}
\left(T > k\,,\,d=i\right)\to 1 \quad\mbox{as}\quad \alpha_\zs{max}\to 0
$$
and using \eqref{sec:01-05-1} and \eqref{sec:Bay.3}, we finally obtain the asymptotic inequality
\[
\EV_\zs{i,k,\theta_\zs{i}}\,\brcs{(T-k)^{r}_\zs{+}\Chi_\zs{\{d=i\}}} \ge
 \brcs{(1-\varepsilon)\frac{|\log\alpha_\zs{j,i}|}{I_\zs{i,j}(\theta_i,\theta_j)}}^r (1+o(1)), \quad \alpha_{max} \to 0
\]
(where $o(1) \to 0$), which holds for an arbitrary $\varepsilon \in(0, 1)$, so  letting $\varepsilon\to 0$ implies lower bounds  \eqref{LW-bnd-jneqi} (for $j\neq i$) and \eqref{LW-bnd-jeqi} (for $j=i$).  The proof is complete.
\endproof

\subsection{The local constraints setting}

To find asymptotic lower bounds for the problems  \eqref{sec:Prbf.5-0}  and \eqref{sec:Prbf.5} in addition to condition $(\A_{1})$ we impose the following condition. 

\noindent $(\H_\zs{1})$ {\em  The parameters $\varrho$, $\m^{*}$ and $\k^{*}$  in \eqref{sec:alph-1-2-16-03-2}  are  functions of $\beta$, i.e. 
$\varrho=\varrho_\zs{\beta}$,
$\m^{*}=\m^{*}_\zs{\beta}$ and $\k^{*}=\k^{*}_\zs{\beta}$,  such that
\begin{equation}\label{sec:Absrsk.4}
\lim_\zs{\beta_\zs{max}\to 0}\varrho_\zs{\beta}=0
\,,\quad
\lim_\zs{\beta_\zs{max}\to 0}\frac{\vert \log \varrho_\zs{\beta}\vert}{\vert \log \alpha^{(2)}_\zs{max}\vert}=0
\quad\mbox{and}\quad
\lim_\zs{\beta_\zs{max}\to 0}\, 
\max_\zs{1\le i,j\le N}\,
\left\vert
\frac{|\log \alpha^{(2)}_\zs{i,j}|}{|\log\beta_\zs{i,j}|}
-
1
\right\vert
=1\,,
\end{equation}
where  $\alpha^{(2)}_\zs{max}=\max_\zs{1\le i,j\le N} \alpha^{(2)}_\zs{i,j}$ and $\beta_\zs{max}=\max_\zs{1\le i,j\le N} \beta_\zs{i,j}$.}

For example, we can take
\begin{equation}\label{sec:Absrsk.9}
\m^{*}_{\beta}=[ \vert\log\beta_\zs{min}\vert/\varrho_\zs{\beta}]
\,,\quad
\k^{*}_{\beta}=\check{\k}\,\m^{*}_\zs{\beta}
\quad\mbox{and}\quad 
\varrho_\zs{\beta}=\frac{1}{1+|\log \beta_\zs{max}|}
\,,
\end{equation}
where $[x]$ is the integer part of the $x$, $\beta_\zs{min}=\min_\zs{1\le i,j\le N}\beta_\zs{i,j}$ and $\check{\k} >1$ is  some fixed number.

The following theorem establishes asymptotic lower bounds in class of detection-identification procedures $\Hc\left(\beta,\k^{*},\m^{*}\right)$.

\begin{theorem} \label{Th.sec:Cnrsk.1} 
 Assume  that  conditions $(\A_{1})$ and $(\H_\zs{1})$ hold. 
Then, for any $r\ge 1$, $k \ge 0$, $1\le i\le N$ and $\theta_i\in\Theta_\zs{i}$, 
\begin{equation} \label{sec:Cnrsk.5}
\liminf_\zs{\beta_\zs{max}\to 0}  \frac{ \inf_\zs{\delta\in \Hc\left(\beta,\k^{*},\m^{*}\right)} 
\,
\Rc_\zs{i,k,\theta_i}(\delta)}{\b^{r}_\zs{i,\beta}(\theta_i)}\,
 \ge 1\,,
\end{equation}
where the denominator $\b_\zs{i,\beta}(\theta_i)$ is defined in \eqref{Lower-bound-1} by replacing the matrix $\alpha$ with $\beta$.
\end{theorem}

\proof
First of all, note that the last condition in
\eqref{sec:Absrsk.4} implies that for any $1\le i\le N$ and $\theta_i\in\Theta_\zs{i}$
$$
\lim_\zs{\beta_\zs{max}\to 0}\frac{\b_\zs{i,\alpha_\zs{2}}(\theta_i)}{\b_\zs{i,\beta}(\theta_i)}
=1\,.
$$
\noindent
It is clear that the last inclusion in \eqref{sec:Comp-16-03} implies that for all $k\ge 0$  and sufficiently small $\beta>0$
 $$
\inf_\zs{\delta\in \Hc\left(\beta,\k^{*},\m^{*}\right)}\,\Rc_\zs{i,k,\theta_i}(\delta)\,\ge 
\inf_\zs{\delta\in \Delta(\alpha_{2},\varrho_\zs{\beta})}\,\Rc_\zs{i,k,\theta_i}(\delta)
\ge\,
\inf_\zs{\delta\in \Delta(\alpha_{2},\varrho_\zs{\beta})}\,
\EV_\zs{i,k,\theta_i} \brcs{(T-k)^{r}_\zs{+}\Chi_\zs{\{d=i\}}}
\,.
$$ 
Now, the lower bounds \eqref{sec:Cnrsk.5} directly  follow from the lower bounds \eqref{sec:Bay.1-nw-1} and  condition $(\H_\zs{1})$. 
\endproof

\section{Upper bounds and asymptotic optimality}\label{sec:Up-bnd}

We begin with studying the sequential procedure \eqref{seq-prs-1} for large threshold values $A_{i,j}$. For any matrix $A=(A_\zs{i,j})_\zs{1\le i,j\le N}$ with $A_\zs{i,j}>1$ and
any  $\vartheta\in\Theta_\zs{i}$ define
\begin{equation}
\label{Up-bnd-21-0-1}
\B_\zs{i,A}(\vartheta)=
\max_\zs{1\le j\le N}\,
\frac{\log A_\zs{i,j}}{\iota_\zs{i,j}(\vartheta)}
\,,
\end{equation}
where the ``information'' functions $\iota_\zs{i,j}(\cdot)$ are defined in \eqref{Lower-bound-1}. We need the following condition:

\noindent $(\H_\zs{2})$ {\em The matrix $A=(A_\zs{i,j})_\zs{1\le i,j\le N}$ is such that 
$$
\lim_\zs{A_\zs{min}\to \infty}
\max_\zs{1\le i\le N}
\left\vert
 \frac{\max_\zs{1\le j\le N}\log A_\zs{i,j}}{\log A_\zs{i,i}}
  -1
\right\vert
=0
\,.
$$
}

\begin{proposition} \label{Prop}
\label{Pr.sec:Up-bnd-1}
If conditions $(\A_\zs{2})$ and $(\H_\zs{2})$ hold true, 
then for any $1\le i\le N$ and any compact set $\K\subset \Theta$ the sequential procedure 
\eqref{seq-prs-1}, in which $\vert\log\varrho\vert=\ao(\log A_\zs{min})$ as $A_\zs{min}\to \infty$, 
admits the following upper bound
\begin{equation}
\label{Up-bnd-21-0-2-NV}
\limsup_\zs{A_\zs{min}\to\infty}\,
 \max_\zs{1\le i\le N}\,  \sup_\zs{\theta_\zs{i}\in \K}\,
 \sup_\zs{1\le k\le \k_\zs{*}}\,
\frac{ \Rc_\zs{i,k,\theta_\zs{i}}(\delta^{*}_\zs{A})}{\B^{r}_\zs{i,A}(\theta_\zs{i})} 
\le 1
\,,
\end{equation}
where
$\theta=(\theta_\zs{1},\ldots,\theta_\zs{N})$,
 $A_\zs{min}=\min_\zs{1\le i,j\le N} A_\zs{i,j}$ and $\k_\zs{*}$ is such that $\k_\zs{*}=\ao(\log A_\zs{min})$ as $A_\zs{min}\to\infty$.
\end{proposition}

\proof
First, note that in view of \eqref{upper-bnd-2}
$\delta^{*}_\zs{A}$ belongs to class $\Delta(\varrho,\alpha)$ with
$$
\alpha_\zs{i,j}=\frac{1}{1+A_\zs{i,i}}\Chi_\zs{\{i=j\}}
+
\frac{1}{1+A_\zs{j,i}}\Chi_\zs{\{i\neq j\}}
\,.
$$
\noindent Therefore, using the upper bound
\eqref{Up-tau-1}, we obtain that uniformly over $1\le k\le \k_\zs{*}$
\begin{align*}
\Pb_\zs{i,k,\theta_\zs{i}}&(T^{*}_\zs{A}>k)=\Pb^{*}(T^{*}_\zs{A}>k)\ge 1-(1-\varrho)^{-k}\sum^{N}_\zs{j=1}\frac{1}{1+A_\zs{j,j}}\\[2mm]
&\ge 1-(1-\varrho)^{-k_\zs{*}}\sum^{N}_\zs{j=1}\frac{1}{1+A_\zs{j,j}}
\to\quad 1\,,\quad A_\zs{min}\to \infty\,.
\end{align*}
Therefore, to obtain  the inequality
\eqref{Up-bnd-21-0-2-NV} it suffices to show that
\begin{equation}
\label{Up-bnd-21-0-2}
\limsup_\zs{A_\zs{min}\to\infty}\,
 \max_\zs{1\le i\le N}\, \sup_\zs{\theta_i\in \K}\,
 \sup_\zs{1\le k\le \k_\zs{*}}\,
\frac{ 
\EV_\zs{i,k,\theta_\zs{i}}\,[(T^{*}_\zs{A}-k)^{r}_\zs{+}]\Chi_\zs{\{d^{*}_\zs{A}=i\}}}{\B^{r}_\zs{i,A}(\theta_i)}
\le 1\,.
\end{equation}

Note also that by condition  $(\A_\zs{2}$) for arbitrary $0<\varepsilon<1$ we can chose such $0<\zeta<1$ for which 
$\left\{\vartheta\,:\,\vert \vartheta-\theta_\zs{i}\vert<\zeta\right\}\subset \Theta_\zs{i}$ for all $1\le i\le N$,
$$
\U^{*}_\zs{1}(\zeta)=
\max_\zs{1\le i\le N}\, \sup_\zs{\theta_i\in\K} 
\sum_{n=1}^\infty \, n^{r-1} \, \sup_\zs{k \ge 0} \Pb_\zs{i,k,\theta_\zs{i}}\brc{
 \inf_\zs{\vert u-\theta_\zs{i}\vert<\zeta} Z_\zs{i,k+n}^{k}(u) < (1  - \varepsilon)\wt{I}_\zs{i}(\theta_\zs{i})n} <\infty
$$
and
$$
\U^{*}_\zs{2}(\zeta)=
\max_\zs{ i\neq j}\, \sup_\zs{\theta_i\in\K} 
\sum_\zs{n=1}^\infty \, n^{r-1} \, \max_\zs{0\le k \le \zeta n} \Pb_\zs{i,k,\theta_\zs{i}}\brc{
 \inf_\zs{\vert u-\theta_\zs{i}\vert<\zeta} Z_\zs{i,k+n}^{k}(u) - Z^{*}_\zs{j,k+n} < (1  - \varepsilon)\wt{I}_\zs{i}(\theta_\zs{i})n} <\infty\,.
$$

Now, for an arbitrary $\vartheta\in\Theta_\zs{i}$, we set
\begin{equation}
\label{sec:Up-bnd-21-3-n}
\n^{*}_\zs{i}
= 
\n^{*}_\zs{i,A}(\vartheta)
=
\left[
\frac{1+\varepsilon}{1-\varepsilon}
\varkappa_\zs{i,A}
\B_\zs{i,A}(\vartheta)
\right]+1
\quad\mbox{and}\quad
\varkappa_\zs{i,A}=\frac{\max_\zs{1\le j\le N} \log A_\zs{i,j}}{\log A_\zs{i,i}}
\,.
\end{equation}
\noindent
Note that in view of the properties
\eqref{sec:Main-Cnds-Lb-1Ps-1}
and the fact that $\k_\zs{*}=\ao(\log A_\zs{min})$ we can conclude that
 for a sufficiently large $A_\zs{min}$ we have $\k_\zs{*}\le \zeta\n^{*}_\zs{i}$. Moreover, 
\begin{equation}
\label{sec:Up-bnd-Mn-1}
\EV_\zs{i,k,\theta_i}\,\brcs{(T^{*}_\zs{A}-k)^{r}_\zs{+}\Chi_\zs{\{d^{*}_\zs{A}=i\}}} \le
\sum_\zs{n\ge 0}\,n^{r-1}\Pb_\zs{i,k,\theta_i}\,\left(T^{*}_\zs{i,A}> k+n\right)
\le
1+(\n^{*}_\zs{i})^{r}
+
\sum_\zs{n> \n^{*}_\zs{i}}n^{r-1}
\Pb_\zs{i,k,\theta_i}\,\left(T^{*}_\zs{i,A}> k+n\right)\,.
\end{equation}
\noindent
Now, the definition \eqref{WSR-def} implies that
$$
\Pb_\zs{i,k,\theta_i}\,\left(T^{*}_\zs{i,A}> k+n\right)
=
\Pb_\zs{i,k,\theta_i}\,
\left(
\max_\zs{1\le l \le n+k}
\min_\zs{1\le j \le N} \frac{<\U_\zs{l}>_\zs{i,j}}{A_\zs{i,j}}
< 1
\right)
\le 
\sum^{N}_\zs{j=1} 
\,
\Pb_\zs{i,k,\theta_i}\,
\left(<\U_\zs{k+n}>_\zs{i,j}
< A_\zs{i,j}
\right)
\,,
$$
i.e. 
\begin{equation}
\label{sec:Up-bnd-Mn*-2}
\Pb_\zs{i,k,\theta_i}\,\left(T^{*}_\zs{i,A}> k+n\right)
\le 
\sum^{N}_\zs{j=1} 
\,
\Pb_\zs{i,k,\theta_i}\,
\left(\log<\U_\zs{k+n}>_\zs{i,j}
<\log A_\zs{i,j}
\right)
\,.
\end{equation}

\noindent
Using \eqref{ST-RdMatr} we obtain that for $i\neq j$
\begin{equation} 
\label{sec:Upp-bnd-F1}
\Pb_\zs{i,k,\theta_i}\,
\left(\log<\U_\zs{k+n}>_\zs{i,j}
<\log A_\zs{i,j}
\right)
=
\Pb_\zs{i,k,\theta_i}\,
\left(\log\L_\zs{i,k+n}-
\log\wh{\L}_\zs{j,k+n}
<\log A_\zs{i,j}
\right)
\end{equation}
and for $i=j$
\begin{equation} 
\label{sec:Upp-bnd-S2}
\Pb_\zs{i,k,\theta_i}\,
\left(\log<\U_\zs{k+n}>_\zs{i,i}
<\log A_\zs{i,i}
\right)
=
\Pb_\zs{i,k,\theta_i}\,
\left(\log\L_\zs{i,k+n}-
(k+n)\log (1-\varrho)
<\log A_\zs{i,i}
\right)\,.
\end{equation}

\noindent
Note that 
$$
\log
\L_\zs{i,k+n}
\ge
\log\pi_\zs{k}\int_\zs{\{\vert  \upsilon -\theta_i\vert<\zeta\}}\,\g^{*}_\zs{i,k,k+n}(\upsilon) W_\zs{i}(\d \upsilon)
\ge \inf_\zs{\{\vert  \upsilon -\theta_i\vert<\zeta\}} Z^{k}_\zs{i,k+n}(\upsilon)
+\l_\zs{i}(\varrho)\,,
$$
where $\l_\zs{i}(\varrho)=-\log \varrho- \k_\zs{*} \log (1-\varrho) - \log W_\zs{i}\left(z\in\Theta_\zs{i}\,:\,\vert z-\theta_i \vert <\zeta\right)$.
Since
$\vert\log\varrho\vert=\ao(\log A_\zs{min})$ and $\k_\zs{*}=\ao(\log A_\zs{min})$, we get 
$$
\lim_\zs{A_\zs{min}\to\infty}\,
\frac{\max_\zs{1\le i\le N}\l_\zs{i}(\varrho)}{\log A_\zs{min}}
=0
\,.
$$
Obviously,
$
\log\wh{\L}_\zs{j,k+n}
\le
Z^{*}_\zs{j,k+n}$. 
Moreover, taking into account that for any $1\le j\le N$
$$
\n^{*}_\zs{i}
\ge \frac{(1+\varepsilon)\log A_\zs{i,j}}{(1-\varepsilon)\wt{I}_\zs{i}(\theta_i)}
$$
 we can obtain that for any $j\neq i$, $\theta_i\in\Theta_\zs{i}$, $n> \n^{*}_\zs{i}>\k_\zs{*}/\zeta$ and 
$0\le k\le \k_\zs{*}$ and for sufficiently large $A_\zs{min}$ for which $\max_\zs{1\le i\le N}\l_\zs{i}(\varrho) \le  \varepsilon \log A_\zs{min}$
\begin{align*}
\Pb_\zs{i,k,\theta_i}\,
&\left(\log<\U_\zs{k+n}>_\zs{i,j}
<\log A_\zs{i,j}
\right)
\le 
\Pb_\zs{i,k,\theta_i}\,
\left(
\inf_\zs{\{\vert  \upsilon -\theta_i\vert<\zeta\}} Z^{k}_\zs{i,k+n}(\upsilon)-
Z^{*}_\zs{j,k+n}
<\log A_\zs{i,j}+\l_\zs{i}(\varrho)
\right)\\[2mm]
&\le 
\max_\zs{0\le l \le \zeta n}\,
\Pb_\zs{i,l,\theta_i}\left(
 \inf_\zs{\vert \upsilon-\theta_i\vert<\zeta} Z_\zs{i,l+n}^{l}(\upsilon) - Z^{*}_\zs{j,l+n} < (1  - \varepsilon)\wt{I}_\zs{i}(\theta_i) n
 \right)
\end{align*}
and
\begin{align*}
\Pb_\zs{i,k,\theta_i}\,&
\left(\log<\U_\zs{k+n}>_\zs{i,i}
<\log A_\zs{i,i}
\right)
\le 
\Pb_\zs{i,k,\theta_i}\,
\left(
\inf_\zs{\{\vert  \upsilon -\theta_i\vert<\zeta\}} Z^{k}_\zs{i,k+n}(\upsilon)
<\log A_\zs{i,i}
+
\l_\zs{i}(\varrho)
\right)\\[2mm]
&\le 
\max_\zs{0\le l \le \zeta n}\,
\Pb_\zs{i,l,\theta_i}\left(
 \inf_\zs{\vert \upsilon-\theta_i\vert<\zeta} Z_\zs{i,l+n}^{l}(\upsilon)  < (1  - \varepsilon)\wt{I}_\zs{i}(\theta_i) n
 \right)\,.
\end{align*}
Therefore, from
\eqref{sec:Up-bnd-Mn-1} we get
\begin{align*}
\EV_\zs{i,k,\theta_i}\,\brcs{(T^{*}_\zs{A}-k)^{r}_\zs{+}\Chi_\zs{\{d^{*}_\zs{A}=i\}}}
&\le
1+(\n^{*}_\zs{i})^{r}
+
\U^{*}_\zs{1}(\zeta)
+\U^{*}_\zs{2}(\zeta)\\[2mm]
&\le 
1+
\left(
\frac{1+\varepsilon}{1-\varepsilon}
\right)^{r}
\varkappa^{r}_\zs{i,A}
\B^{r}_\zs{i,A}(\theta_i)
+
\U^{*}_\zs{1}(\zeta)
+\U^{*}_\zs{2}(\zeta)
\,,
\end{align*}
and using the condition $(\H_\zs{2})$, we get
$$
\limsup_\zs{A_\zs{min}\to\infty}\,
 \max_\zs{1\le i\le N}\,  \sup_\zs{\theta_i\in \K}\,
 \sup_\zs{1\le k\le \k_\zs{*}}\,
\frac{ 
\EV_\zs{i,k,\theta_\zs{i}}\,\brcs{(T^{*}_\zs{A}-k)^{r}_\zs{+}\Chi_\zs{\{d^{*}_\zs{A}=i\}}}}{\B^{r}_\zs{i,A}(\theta_i)}
\le 
\left(
\frac{1+\varepsilon}{1-\varepsilon}
\right)^{r}.
$$
Since $\varepsilon$ can be arbitrarily small, taking the limit as $\varepsilon \to 0$, we obtain the bound \eqref{Up-bnd-21-0-2}, which completes the proof of 
Proposition \ref{Pr.sec:Up-bnd-1}.
\endproof

\begin{remark}
If both left-tail and right-tail conditions $(\A_\zs{1})$ and  $(\A_\zs{2})$  hold along with conditions $(\H_\zs{1})$ and $(\H_\zs{2})$, then
inverting the equality \eqref{upper-bnd-2-1**} and using Theorem \ref{Th.sec:Cnrsk.1} (with $\beta$ replaced with $A^{-1}$) and Proposition~\ref{Prop} simultaneously it can be shown that the following 
asymptotic equalities for the moments of delay of the procedure $\delta_A$ hold for any fixed $k$, $\theta_i\in \Theta_i$ and all $i=1,\dots,N$:
\begin{equation*}
\Rc_{i,k,\theta_i}(\delta_A) = \max_\zs{1\le j\le N}\,
\frac{\log A_\zs{i,j}}{\iota_\zs{i,j}(\theta_i)} (1+o(1)) = \max\set{ \frac{\log A_\zs{i,i}}{\widetilde{I}_\zs{i}(\theta_i)},  \max_\zs{1\le j \neq i\le N}\,
\frac{\log A_\zs{i,j}}{\widehat{I}_\zs{i,j}(\theta_i)}} (1+o(1))  \quad \text{as $A_{min}\to \infty$}.
\end{equation*}
\end{remark}

To obtain the optimal detection rate we need to impose the following condition:

\noindent $(\H_\zs{3})$  {\em 
Parameters $\varrho^{opt}$ and $\k^{*}$ are functions of $\beta$, i.e.
$\varrho^{opt}=\varrho^{opt}_\zs{\beta}$,
 $\k^{*}=\k^{*}_{\beta}$ and $\m^{*}=\m^{*}_{\beta}$, such that
\begin{equation}\label{sec:Cnrsk.6}
\lim_\zs{\beta_\zs{max}\to 0}\,
\frac{\vert\log\varrho^{opt}\vert}{\vert\log \beta_\zs{max} \vert}
=0
\quad\mbox{and}\quad
\lim_\zs{\beta_\zs{max}\to 0}\, 
\max_\zs{1\le i,j\le N}\,
\left\vert
\frac{|\log \alpha^{(1)}_\zs{i,j}|}{|\log\beta_\zs{i,j}|}
-
1
\right\vert
=0\,,
\end{equation}
where 
$$
\alpha^{(1)}_\zs{i,j}
=\beta_\zs{i,i}(1-\varrho^{opt})^{\k^{*}}\Chi_\zs{\{i= j\}}
+
\beta_\zs{i,j}\varrho^{opt}(1-\varrho^{opt})^{\k^{*}}\Chi_\zs{\{i\neq j\}}
\,.
$$
}

\noindent
For example, for some $\check{\k}>1$ we can take
\begin{equation}\label{sec:Up-Bnd.1}
\m^{*}_\zs{\beta}=\left[\frac{ \vert\log\beta_\zs{min}\vert}{\varrho_\zs{\beta}}
\right]
\,,\quad
\k^{*}_{\beta}=\check{\k}\m^{*}_\zs{\beta}\,,\quad
\varrho^{opt}
=  \frac{\vert \log \beta_\zs{max}\vert \varrho_\zs{\beta}}{\vert \log \beta_\zs{min}\vert(1+\vert\log \varrho_\zs{\beta} \vert)}
\quad\mbox{and}\quad 
\varrho_\zs{\beta}=\frac{1}{1+|\log \beta_\zs{max}|}
\,.
\end{equation}
Then under the conditions
\begin{equation}
\label{sec:Up-Bn-121-21}
\lim_\zs{\beta_\zs{max}\to 0}\frac{\log \vert\log\beta_\zs{min}\vert}{\vert\log\beta_\zs{max}\vert}
=0
\quad\mbox{and}\quad
\lim_\zs{\beta_\zs{max}\to 0}\,
\max_\zs{1\le i\le N}
\left\vert
\frac{\max_\zs{1\le j\le N}\vert\log\beta_\zs{j,i} \vert}{\vert \log\beta_\zs{i,i}\vert}
-
1
\right\vert=0
\,,
\end{equation}
we obtain that the conditions $(\H_\zs{1})$ -- $(\H_\zs{3})$ hold.

Denote by $\delta^{opt}_\zs{\beta}=(T^{opt},d^{opt})$   the procedure \eqref{seq-prs-1} with $\varrho=\varrho^{opt}$ and
\begin{equation}\label{sec:Cnrsk.7}
 A_\zs{i,j}= A_\zs{i,j}^{opt}(\beta)
=
\left(
\frac{1+\tr \beta}{\beta_\zs{i,i}(1-\varrho^{opt})^{\k^{*}}}
-1
\right)
\Chi_\zs{\{i= j\}}
+
\frac{1+\tr \beta}{\beta_\zs{j,i}
 (\varrho^{opt}) (1-\varrho^{opt})^{\k^{*}}}
 \Chi_\zs{\{i\neq j\}} \, ,
\,
\end{equation}
that is,
\begin{equation}
\label{seq-prs-opt}
T^{opt}_\zs{\beta}=\min_\zs{1\le i\le N}\,T^{opt}_\zs{i,\beta} \quad\mbox{and}\quad d^{opt}_\zs{\beta}=i
\quad\mbox{if}\quad T^{opt}_\zs{i,\beta}=T^{opt}_\zs{\beta}\,,
\end{equation}
where
\begin{equation}\label{Topt}
T^{opt}_\zs{i,\beta}=\inf\set{n \ge 1: \min_\zs{1\le j\le N}\, \frac{<\U_\zs{n}>_\zs{i,j}}{A_\zs{i,j}^{opt}(\beta)}\ge 1}
 \,.
\end{equation}

The following theorem deduces the pointwise and minimax optimality properties of  the procedure $\delta^{opt}_\zs{\beta}$.

\begin{theorem}\label{Th.sec:Cnrsk.2} 
 Assume that conditions $(\A_\zs{1})$--$(\A_\zs{2})$ and $(\H_\zs{1})$--$(\H_\zs{3})$ hold true. Then 
the procedure $\delta^{opt}_\beta$ is optimal in the pointwise sense, i.e., for any $\theta_i\in\Theta_\zs{i}$, $1\le i\le N$,
 and for every  fixed $k\ge 0$
\begin{equation}\label{sec:SRAO2}
\lim_\zs{\beta_\zs{max}\to 0}
\frac{ \inf_\zs{\delta\in \Hc\left(\beta,\k^{*},\m^{*}\right)}\, \Rc_\zs{i,k,\theta_i}(\delta)}{\Rc_\zs{i,k,\theta_i}(\delta^{opt}_\zs{\beta})}
=1
\quad\mbox{and}\quad
\lim_\zs{\beta_\zs{max}\to 0}
\frac{\Rc_\zs{i,k,\theta_i}(\delta^{opt}_\zs{\beta})}{\b^{r}_\zs{i,\beta}(\theta_i)}=1\,.
\end{equation}
Also, for any $\k_\zs{*}=\ao(\vert\log\beta_\zs{max}\vert)$ as $\beta_\zs{max}\to 0$
 the procedure $\delta^{opt}_\zs{\beta}$ is optimal 
in the minimax sense, i.e.,  for any  $\theta_i\in\Theta_\zs{i}$ and $1\le i\le N$,
\begin{equation}\label{sec:SRAO2b}
\lim_\zs{\beta_\zs{max}\to 0}
\frac{ \inf_\zs{\delta\in \Hc\left(\beta,\k^{*},\m^{*}\right)}\, \max_\zs{1\le k\le \k_\zs{*}}\Rc_\zs{i,k,\theta_i}(\delta)}
{ \max_\zs{1\le k\le \k_\zs{*}}\Rc_\zs{i,k,\theta_i}(\delta^{opt}_\zs{\beta})}
=1
\quad\mbox{and}\quad
\lim_\zs{\beta_\zs{max}\to 0}
\frac{ \max_\zs{1\le k\le \k_\zs{*}}\Rc_\zs{i,k,\theta_i}(\delta^{opt}_\zs{\beta})}{\b^{r}_\zs{i,\beta}(\theta_i)}=1\,.
\end{equation}
\end{theorem}

\proof
By condition $(\H_\zs{3})$ 
$$
\lim_\zs{\beta_\zs{max}\to 0}\, 
\max_\zs{1\le i,j\le N}\,
\left\vert
\frac{|\log A_\zs{i,j}^{opt}(\beta)|}{|\log\beta_\zs{j,i}|}
-
1
\right\vert
=0,
$$
so that using the asymptotic upper bound \eqref{Up-bnd-21-0-2-NV} in Proposition~\ref{Prop} we obtain the asymptotic upper bound  
\begin{equation}
\label{Up-bnd-21-0-2-NVbeta}
\limsup_\zs{\beta_\zs{max}\to 0}\,
 \max_\zs{1\le i\le N}\,  \sup_\zs{\theta_i\in \K}\,
 \sup_\zs{1\le k\le \k_\zs{*}}\,
\frac{ \Rc_\zs{i,k,\theta_\zs{i}}(\delta_\beta^{opt})}{\B^{r}_\zs{i,\beta}(\theta_i)} 
\le 1
\,.
\end{equation}
Comparing this bound with the lower bound \eqref{sec:Bay.1-nw-1} in Theorem~\ref{Th.sec:Bay.1} yields \eqref{sec:SRAO2} and \eqref{sec:SRAO2b}.
\endproof

The next theorem also shows that the procedure $\delta^{opt}_\zs{\beta}=(T^{opt},d^{opt})$ is ``robust'' in the following sense
\begin{equation}
\label{Robust-risk-*}
\Rc^{*}_\zs{\beta}(\delta)=
\sup_\zs{\theta\in\Theta}\,
\max_\zs{1\le i\le N}\,
\max_\zs{1\le k\le \k_\zs{*}}\,
\frac{\Rc_\zs{i,k,\theta_i}(\delta)}{\b^{r}_\zs{i,\beta}(\theta_i)}
\,. 
\end{equation}

\begin{theorem}\label{Th.sec:Cnrsk.2-*} 
 Suppose that  conditions $(\A_\zs{1})$--$(\A_\zs{2})$  and $(\H_\zs{1})$--$(\H_\zs{3})$ hold 
 and  $\k_\zs{*}=\ao(\vert\log\beta_\zs{max}\vert)$ as $\beta_\zs{max}\to 0$.
 Then 
\begin{equation}\label{sec:SRAO2c}
\lim_\zs{\beta_\zs{max}\to 0}
\frac{ \inf_\zs{\delta\in \Hc\left(\beta,\k^{*},\m^{*}\right)}\, \Rc^{*}_\zs{\beta}(\delta)}{\Rc^{*}_\zs{\beta}(\delta^{opt}_\zs{\beta})}
=1
\,.
\end{equation}
\end{theorem}

The proof is similar to the proof of Theorem~2 in \cite{PerTar-JMVA2019} for the single-stream detection problem and is omitted.

%
\section{Detection-identification of changes in homogeneous Markov models}
\label{sec:Mrk}

Let the observations 
 $(X_\zs{i,n})_\zs{n\ge 1}$ be time homogeneous  Markov processes with values in a measurable space $(\Xc_\zs{i},\Bc_\zs{i})$
defined by a family
of the transition probabilities
$(P^{\theta_i}_\zs{i}(x,A))_\zs{\theta_i\in\Theta_\zs{i}}$ for some fixed parameter set $\Theta_\zs{i}\subseteq \bbr^{p}$.  
In the sequel we denote by $\EV^{\theta_i}_\zs{i,x}(\cdot)$   the  expectation 
with respect to this probability.
Moreover, we assume that for any $1\le i\le N$ the observations
 $(X_\zs{i,n})_\zs{n\ge 1}$ are  Markov processes, such that
$(X_\zs{i,n})_\zs{1\le n \le \nu}$ is a homogeneous process
 with the transition (from $x$ to $y$) density $f^{*}_\zs{i}(y|x)$ and in the case when $\nu=+\infty$
 this process is ergodic with the ergodic distribution $\lambda_\zs{*,i}$. 
 We denote by $\Pb^{*}$
 the distribution 
of the observations $(X_\zs{i,n})_\zs{1\le i\le N,n\ge 1}$ 
 of this process when $\nu=\infty$. The expectation with respect to this distribution will be denoted by $\EV^{*}(\cdot)$.
In addition, we assume that for any $1\le i\le N$ the process $(X_\zs{i,n})_\zs{n> \nu}$ is 
homogeneous positive  ergodic  with the transition density $f_\zs{i,\theta_i}(y|x)$ 
and the ergodic (stationary) distribution $\lambda_\zs{\theta_i,i}$ ($\theta_i\in\Theta_\zs{i}$). 
The densities $f^{*}_\zs{i}(y|x)$ and  $(f_\zs{i,\theta_i}(y|x))_\zs{\theta_\in\Theta_\zs{i}}$   are calculated 
with respect to a sigma-finite positive measure $\mu_\zs{i}$ on $\Bc_\zs{i}$.
In this case, we can represent the LLR process $Z^{k}_\zs{i,n}(u)$ defined in \eqref{Znk_df} as
\begin{equation}\label{sec:Mrk.5}
Z^{k}_\zs{i,n}(u)=\sum^{n}_\zs{j=k+1} g_\zs{i}(u,X_\zs{i,j},X_\zs{i,j-1})\,,
\quad
  g_\zs{i}(u,y,x)=\log \frac{f_\zs{i,u}(y|x)}{f^{*}_\zs{i}(y|x)}\,.
\end{equation}
\noindent

We also assume that densities $f_\zs{i,u}(y|x)$ are continuously differentiable 
with respect to $u$ in a compact set  $\widetilde{\Theta}_\zs{i}\subseteq\Theta_\zs{i}$.  
Now we set
\begin{equation}\label{sec:Mrk.6-h-def}
h(x,y)=\max_\zs{1\le i\le N}
\max_\zs{u\in\widetilde{\Theta}_\zs{i}}\,
\max_\zs{1\le j\le p}\,
\vert \partial g_\zs{i}(u,y,x)/
\partial u_\zs{j}
\vert 
\end{equation}
and
$$
\overline{h}_\zs{i}(\theta_i,y)=
\int_\zs{\Xc_\zs{i}}\,h(y,x)\,f_\zs{i,\theta_i}(y|x)\,\mu_\zs{i}(\d y)
\,,\quad \theta_i\in\Theta_\zs{i}\,.
$$
\noindent
 For some $\q>0$ define
$$
g^{*}_\zs{\q}(x)=
\sup_\zs{n\ge 1}\,
\max_\zs{1\le i\le N}\,
\sup_\zs{\theta_i\in\Theta_\zs{i}}\,
\EV^{\theta_i}_\zs{i,x}\,|g_\zs{i}(\theta,X_\zs{i,n},X_\zs{i,n-1})|^{\q}
$$
and
\begin{equation}
\label{sec:Mrk.2}
h^{*}_\zs{\q}(x)=\sup_\zs{n\ge 1}\,
\max_\zs{1\le i\le N}\,
\sup_\zs{\theta_i\in\Theta_\zs{i}}\,\EV^{\theta_i}_\zs{i,x}\,|h(X_\zs{i,n},X_\zs{i,n-1})|^{\q}
\,.
\end{equation}
Also, define
\begin{equation}\label{sec:Mrk.6}
J_\zs{i}(\theta_i,x)=
\int_\zs{\Xc_\zs{i}}\,g_\zs{i}(\theta_i,y,x)\,f_\zs{i,\theta_i}(y|x)\,\mu_\zs{i}(\d y)
\quad\mbox{and}\quad
J^{*}_\zs{i}(\theta_i,x)=
\int_\zs{\Xc_\zs{i}}\,g_\zs{i}(\theta_i,y,x)\,f^{*}_\zs{i}(y|x)\,\mu_\zs{i}(\d y) 
\,.
\end{equation}
Obviously, $J_\zs{i}(\theta_i,x)\ge 0$ and $J^{*}_\zs{i}(\theta_i,x)\le 0$. Write
\begin{equation}\label{sec:Mrk.6-erg}
\overline{J}_\zs{i}(\theta_i)
=
\int_\zs{\Xc_\zs{i}}\,J_\zs{i}(\theta_i,x) \lambda_\zs{\theta_i,i}(\d x)
\quad\mbox{and}\quad
\overline{J}^{*}_\zs{i}(\theta_i)
=
\int_\zs{\Xc_\zs{i}}\,J^{*}_\zs{i}(\theta_i,x) \lambda_\zs{*,i}(\d x) 
\,.
\end{equation}

Introduce the following conditions.

\noindent 
$(\C_\zs{1})$ {\em For any $1\le i\le N$ there exist  sets $C_\zs{i}\in\Bc_\zs{i}$ with $\mu_\zs{i}(C_\zs{i})<\infty$ such that
\begin{enumerate}
 
 \item[$({\rm C}1.1)$] 
 $
 f_\zs{*}=
 \min_\zs{1\le i\le N}\,
 \inf_\zs{\theta_i\in\Theta_\zs{i}}\,
 \inf_\zs{x,y\in C_\zs{i}}\,f_\zs{i,\theta_i}(y|x)>0$.
\item[$({\rm C}1.2)$] 
For any $1\le i\le N$ there exists  $\Xc_\zs{i}\to [1,\infty)$ Lyapunov's function $\V_\zs{i}$ such that 
\begin{itemize}
\item 
$\V_\zs{i}(x)\ge J_\zs{i}(\theta_i,x)$ and
$\V_\zs{i}(x)\ge \overline{h}_\zs{i}(\theta_i,x)$ for any $\theta_i\in\Theta_\zs{i}$ and $x\in\Xc_\zs{i}$ .

\item
 $
 \max_\zs{1\le i\le N}\,
 \sup_\zs{x\in C_\zs{i}} \V_\zs{i}(x)<\infty$.

\item  There exist $0<\rho<1$ and $\m^{*}>0$ 
such that  for all $1\le i\le N$,  $x\in\Xc_\zs{i}$ and $\theta_i\in\Theta_\zs{i}$, 
\begin{equation}
\label{drift_ineq_11}
\EV^{\theta_i}_\zs{i, x}[\V_\zs{i}(X_\zs{i,1})]
\le 
(1-\rho) \V_\zs{i}(x)
+ \m^{*}\,
\Chi_\zs{C_\zs{i}}(x)\,.
\end{equation}
\end{itemize}
\end{enumerate}
\noindent 
$(\C_\zs{2}(\q))$ 
There exists $\q> 2$ such that for any $1\le i\le N$
$$
\sup_\zs{k\ge 1}\,
\EV^{*}\,[g^{*}_\zs{\q}(X_\zs{i,k})]\,
<\infty
\,,\quad
\sup_\zs{k\ge 1}\,
\EV^{*} [h^{*}_\zs{\q}(X_\zs{i,k})]\,
<\infty
\quad\mbox{and}\quad
\sup_\zs{k\ge 1}\,
\EV^{*}[\upsilon^{*}_\zs{\q}(X_\zs{i,k})]\,
<\infty
\,,
$$
where the functions
 $g^{*}_\zs{\q}(x)$ and
  $h^{*}_\zs{\q}(x)$
  are given in \eqref{sec:Mrk.2}
 and
 \begin{equation}\label{sec:Mrk.3}
\upsilon^{*}_\zs{\q}(x)=\sup_\zs{1\le i\le N\,,\,n\ge 0}\,
\sup_\zs{\theta_i\in\Theta_\zs{i}}\,
\EV^{\theta_i}_\zs{i,x}\,\left[\V_\zs{i}(X_\zs{i,n})\right]^{\q}\,.
\end{equation}
\noindent 
$(\C_\zs{3}(\q))$ The function $g_\zs{i}(u,y,x)$ can be represented as
\begin{equation}
\label{funtc-g}
g_\zs{i}(u,y,x)
-\overline{J}^{*}_\zs{i}(u)
=\sum^{m}_\zs{l=1}\alpha_\zs{i,l}(u)\check{g}_\zs{i,l}(y,x)
\end{equation}
with $ \alpha_\zs{j,l}(u)$ and $\check{g}_\zs{i,l}(y,x)$ such that for any $1\le i\le N$ and $1\le l\le m$
$$
\sup_\zs{u\in\Theta_\zs{i}}
\vert \alpha_\zs{j,l}(u) \vert
<\infty
\quad\mbox{and}\quad
\sup_\zs{n\ge 1} n^{-\q/2}
\EV^{*} \left\vert \sum^{n}_\zs{j=1} \check{g}_\zs{i,l}(X_\zs{i,j},X_\zs{i,j-1})
\right\vert^{\q}
<\infty\,.
$$
}

\begin{theorem} \label{Th.sec:Mrk.1} 
Assume that conditions $(\C_\zs{1})-(\C_\zs{3}(\q))$ hold true and the functions $\overline{J}_\zs{i}(\theta_i)$ and
$-\overline{J}^{*}_\zs{i}(\theta_i)$ defined in \eqref{sec:Mrk.6} are continuous and positive for $\theta_i\in\Theta_\zs{i}$.  
 Then  conditions $(\A_\zs{1})$ and  $(\A_\zs{2}(r))$ are satisfied for any $0<r<\q/2$ with 
\begin{equation}
\label{KL-Inf-1}
I_\zs{i,j}(\theta_\zs{i},\theta_\zs{j})
=
\overline{J}_\zs{i}(\theta_\zs{i})\,
\Chi_\zs{\{i=j\}}
+
[
\overline{J}_\zs{i}(\theta_\zs{i})
-
\overline{J}^{*}_\zs{j}(\theta_\zs{j})
]
\Chi_\zs{\{i\neq j\}}
\,.
\end{equation}
\end{theorem}

\proof
Note first that conditions \eqref{sec:Cnd-1} and  \eqref{rcompLeft-1} follow from Theorem 8 in \cite{PerTar-JMVA2019} 
that uses the uniform geometric ergodicity property and concentration inequalities methods developed
in \cite{GaPe13, GaPe14}. To prove condition \eqref{rcompLeft-2} we observe that condition $(\C_\zs{3})$ and the Chebyshev inequality imply that for any $\a>0$
\begin{equation}
\label{upper-bnd-NAS-1}
\sup_\zs{n\ge 1}\,n^{\q/2}
\max_\zs{1\le i\le N}\,
\max_\zs{1\le k\le n}
\Pb^{*}\brc{\sup_\zs{\theta_i\in\Theta_\zs{i}} \vert \Delta^{k}_\zs{i,n}(\theta_i)\vert > \a n}
<\infty
\,,
\end{equation}
where 
$$
\Delta^{k}_\zs{i,n}(\theta_i) 
=Z^{k}_\zs{i,n} - \overline{J}^{*}_\zs{i}(\theta_i)(n-k)=
\sum^{n}_\zs{j=k+1}\left( g_\zs{i}(\theta_i,X_\zs{i,j},X_\zs{i,j-1})
-\overline{J}^{*}_\zs{i}(\theta_i)
\right)
\,.
$$
Taking into account  that $\overline{J}^{*}_\zs{j}(\theta_i)\le 0$, we obtain that for any $\varepsilon>0$ and $1\le j\le N$
$$
\Pb^{*}\left( Z^{*}_\zs{j,n}>\varepsilon n \right)
\le 
\Pb^{*}\brc{\max_\zs{1\le k\le n}\sup_\zs{\theta_j\in\Theta_\zs{j}} \vert \Delta^{k}_\zs{j,n}(\theta_j)\vert > \varepsilon n}
\le 
n
\max_\zs{1\le k\le n}
\Pb^{*}\brc{\sup_\zs{\theta_j\in\Theta_\zs{j}} \vert \Delta^{k}_\zs{j,n}(\theta_j)\vert > \varepsilon n}
\,.
$$
Therefore, \eqref{upper-bnd-NAS-1} implies the condition \eqref{rcompLeft-2} for $r<\q/2$.
\endproof

\begin{remark}
\label{Re.sec:Mrk-1-0}
The function $\overline{J}_\zs{i}(\cdot)$ is called the
Kullback-Leibler divergence for the Markov processes (see, e.g., \cite{GirardinKonevPergamenshchikov2018}).
\end{remark}

Note that  condition $({\rm C}1.1)$ does not always hold for the process $(X_\zs{i,n})_\zs{n\ge 1}$ directly.
For example, this condition does not hold for the practically important autoregression process of the order more than one.
For this reason, we need to weaken this requirement. Similarly to \cite{PerTar-SISP2016} we assume that there exists $p\ge 2$
for which the  process $(\wt{X}^{\iota}_\zs{i,n})_\zs{n\ge \wt{\nu}}$ 
for $\wt{\nu}_\zs{\iota}=\nu/p-\iota$ and $0\le \iota\le p-1$ defined as $\wt{X}^{\iota}_\zs{i,n}=X_\zs{i,n p+\iota}$ satisfies the following conditions:

\noindent 
$(\C^{\prime}_1)$ {\em There exist sets $C_\zs{i}\in\Bc_\zs{i}$ with $\mu_\zs{i}(C_\zs{i})<\infty$ such that

\begin{enumerate}
 
 \item[$(\C^{\prime}1.1)$] 
 $
 \min_\zs{1\le i\le N}\,
 \inf_\zs{1\le \iota\le p}
 \,
 \inf_\zs{\theta_i\in\Theta_\zs{i}}\,
 \inf_\zs{x,y\in C}\,\wt{f}^{\iota}_\zs{i,\theta_i}(y|x)>0\,,
 $
 where $\wt{f}^{\iota}_\zs{i,\theta_i}(y|x)$ is the transition density for the  process
 $(\wt{X}^{\iota}_\zs{i,n})_\zs{n\ge \wt{\nu}_\zs{i}}$.

\item[$(\C^{\prime}1.2)$] 
For any $1\le i\le N$ there exists  $\Xc_\zs{i}\to [1,\infty)$ Lyapunov's function $\V_\zs{i}$ such that 
\begin{equation*}
\max_\zs{1\le i\le N}\,
\max_\zs{1\le j\le p}\,
 \sup_\zs{\theta_i\in\Theta_\zs{i}}\,
\sup_\zs{x\in\Xc_\zs{i}}\,
 \frac{\EV^{\theta_i}_\zs{i,x}\,\brcs{\V_\zs{i}(X_\zs{i,j})}}{\V_\zs{i}(x)}
 \,<\,\infty
\quad\mbox{and}\quad
\sup_\zs{\theta_i\in\Theta_\zs{i}}\,\lambda_\zs{\theta_i,i}(\V_\zs{i})
\,<\,\infty\,.
\end{equation*}

\begin{itemize}
\item 
$\V_\zs{i}(x)\ge J_\zs{i}(\theta_i,x)$ and $\V_\zs{i}(x)\ge \overline{h}_\zs{i}(\theta_i,x)$ for 
$\theta_i\in\Theta_\zs{i}$ and $x\in\Xc_\zs{i}$ and 
$
\max_\zs{1\le i\le N}\,
 \sup_\zs{\theta_i\in\Theta_\zs{i}}\,
 \sup_\zs{x\in C_\zs{i}} \V_\zs{i}(x)<\infty. 
$

\item For any $1\le i\le N$, there exist $0<\rho<1$ and $\m^{*}>0$
such that 
 for all  $x\in\Xc_\zs{i}$, $\theta_i\in\Theta_\zs{i}$, and
$0\le \iota\le p-1$
\begin{equation}\label{drift_in_+++}
\EV^{\theta_i}_{i,x}\left[\V_\zs{i}(\wt{X}^{\iota}_\zs{i,1})\right]
\le  (1-\rho) \V_\zs{i}(x) + \m^{*}\Chi_\zs{C_\zs{i}}(x)
\,.
\end{equation}
\end{itemize}
\end{enumerate}
}

Similarly to Theorem \ref{Th.sec:Mrk.1} we can prove, using  Theorem  9 in \cite{PerTar-JMVA2019}, the following result.

\begin{theorem} \label{Th.sec:Mrk.1-1} 
Assume that conditions $(\C^{\prime}_\zs{1})$ and $(\C_\zs{2}(\q))$ hold,  the processes $(\wt{X}^{\iota}_\zs{i,n})_\zs{n\ge \wt{\nu}}$ satisfy the condition $(\C_\zs{3}(\q))$,  
and
the functions $\overline{J}_\zs{i}(\theta_i)$ and
$-\overline{J}^{*}_\zs{i}(\theta_i)$ defined in \eqref{sec:Mrk.6} are continuous and positive for $\theta_i\in\Theta_\zs{i}$.  
 Then  conditions $(\A_\zs{1})$ and  $(\A_\zs{2}(r))$ hold for any $0<r<\q/2$ with 
 the functions $I_\zs{i,j}$ defined in \eqref{KL-Inf-1}.
\end{theorem}

To check the condition $\C_\zs{3}(\q))$ we need to obtain concentration inequality for the 
 homogenous Markov process  $(X_\zs{i,n})_\zs{ n \ge 1}$ with the transition density $f^{*}_\zs{i}(y|x)$. The following condition is sufficient for this purpose:

\noindent 
$(\C^{*}_\zs{3}(\q))$ {\em For any $1\le i\le N$ there exist sets $C_\zs{i}\in\Bc_\zs{i}$ with $\mu_\zs{i}(C_\zs{i})<\infty$ such that

\begin{enumerate}
 
 \item  
 $
 \min_\zs{1\le i\le N}\,
 \inf_\zs{x,y\in C_\zs{i}}\,f^{*}_\zs{i}(y|x)>0$.
\item 
For any $1\le i\le N$
there exists  $\Xc_\zs{i}\to [1,\infty)$ Lyapunov's function $\V_\zs{i}$ such that

\begin{itemize}

\item
 $
 \max_\zs{1\le i\le N}\,
 \sup_\zs{x\in C_\zs{i}} \V_\zs{i}(x)<\infty$.

\item  There exist $0<\rho<1$ and $\m^{*}>0$ 
such that  for all $1\le i\le N$,  $x\in\Xc_\zs{i}$ and $\theta_i\in\Theta_\zs{i}$, 
\begin{equation}
\label{drift-ineq-11-*}
\EV^{*}_\zs{x}[\V_\zs{i}(X_\zs{i,1})]
\le 
(1-\rho) \V_\zs{i}(x)
+ \\m^{*}\,
\Chi_\zs{C_\zs{i}}(x)\,.
\end{equation}
\item There exists $\q>0$ such that for any $x_\zs{i}\in \Xc_\zs{i}$
\begin{equation}
\label{UB-1-22}
\max_\zs{1\le i\le N}\,
\sup_\zs{j\ge 1}\,
\EV^{*}_\zs{x_\zs{i}}\,[\V^{\q}_\zs{i}(X_\zs{i,j})]
<\infty\,.
\end{equation}

\end{itemize}
\item
The functions $\check{g}_\zs{i,l}(y,x)$ in  \eqref{funtc-g}
are such that for any $1\le i\le N$ and $x\in \Xc_\zs{i}$
$$
\max_\zs{1\le l\le m}
\left| \int_\zs{\Xc_\zs{i}} \check{g}_\zs{i,l}(y,x)f^{*}_\zs{i}(y|x)\mu_\zs{i}(\d y) 
\right|
\le
\V_\zs{i}(x)
\quad\mbox{and}\quad
\sup_\zs{j\ge 1}\,\EV^{*}_\zs{x} |\check{g}_\zs{i,l}(X_\zs{i,j},X_\zs{i,j-1})|^{q}
<\infty
\,.
$$
\end{enumerate}
}

Proposition 1 from \cite{PerTar-SISP2016} provides the following result.   
\begin{proposition}\label{Pr.sec:Mrk.1-1}
The condition 
$(\C^{*}_\zs{3}(\q))$ implies the condition $(\C_\zs{3}(\q))$. 
\end{proposition}

\section{Examples}\label{sec:Ex}

\subsection{Example 1:  Change in the parameters of  multivariate linear difference equations}\label{ssec:Ex1}

Consider the multivariate models in $\bbr^{p}$ given by
\begin{equation}\label{sec:Ex.6-0n}
X_\zs{i,n} = \left(\Gamma^{*}_\zs{i,n}\Ind{n\le \nu}+\Gamma_\zs{i,n}\Ind{n>\nu}\right)\,X_\zs{i,n-1}+w_\zs{i,n}\,,
\end{equation}
where $\Gamma^{*}_\zs{i,n}$ and $\Gamma_\zs{i,n}$ are $p\times p$ random matrixes and $(w_\zs{i,n})_\zs{n\ge 1}$ is an i.i.d. 
sequence of Gaussian random vectors $\Nc(0, Q^{*}_\zs{i})$ in
$\bbr^{p}$ with the positive definite $p\times p$ matrix $Q^{*}_\zs{i}$. 
Assume also that 
\begin{equation}\label{theta-mtrx-10-05-1}
\Gamma^{*}_\zs{i,n}=\theta^{*}_\zs{i}+B_\zs{i,n}
\quad\mbox{and}\quad
\Gamma_\zs{i,n}=\theta_\zs{i}+B_\zs{i,n}
\end{equation}
 and $(B_\zs{i,n})_\zs{n\ge 1}$ are i.i.d. Gaussian random matrices $\Nc(0\,,Q_\zs{i})$,
where the $p^{2}\times p^{2}$ matrix $Q_\zs{i}=\E[B_\zs{i,1}\otimes B_\zs{i,1}]$ is  positive definited. 
Assume, in addition, that all eigenvalues of
the matrix 
$$
\EV[\Gamma^{*}_\zs{i,1}\,\otimes\,\Gamma^{*}_\zs{i,1}] =\theta^{*}_\zs{i}\otimes \theta^{*}_\zs{i}+Q_\zs{i}
$$
are less than one in module. Define
\begin{equation}\label{theta-set-mtrx}
\Theta^{st}_\zs{i}=\{\theta_i\in\bbr^{p^{2}}\,:\, \max_\zs{1\le j\le p^{4}_\zs{i}} \e_\zs{j}(\theta_i\otimes\theta_i+Q_\zs{i})<1\}\}
\setminus \,
\{\theta^{*}_\zs{i}\}
\,,
\end{equation}
where $\e_\zs{j}(\Gamma)$ is the $j$th eigenvalue of matrix $\Gamma$, and assume further that 
in
\eqref{theta-mtrx-10-05-1}
the matrices
 $\theta_\zs{i}\in\Theta^{st}_\zs{i}$.
In this case,  the processes  $(X_\zs{i,n})_\zs{n\ge 1}$ (for $\nu=\infty$) and $(X_\zs{i,n})_\zs{n> \nu}$ (for $\nu<\infty$)
are ergodic with the ergodic distributions given by the vectors [\citep{KlPe04}]
$$
\varsigma^{*}_\zs{i}=\sum_\zs{l\ge 1}\prod^{l-1}_\zs{j=1}\Gamma^{*}_\zs{i,j}\,w_\zs{i,l}
\quad\mbox{and}\quad
\varsigma_\zs{i,\theta}=\sum_\zs{l\ge 1}\prod^{l-1}_\zs{j=1}\Gamma_\zs{i,j}\,w_\zs{i,l}
$$
i.e.,  the corresponding invariant measures $\lambda_\zs{*,i}$ and $\lambda_\zs{\theta,i}$ on $\bbr^{p}$ are 
defined as 
$$
\lambda_\zs{*,i}(\d x)=\Pb(\varsigma^{*}_\zs{i}\in \d x)
\quad\mbox{and}\quad
\lambda_\zs{\theta,i}(\d x)=\Pb(\varsigma_\zs{i,\theta}\in \d x)
$$
Note that in this case the Markov processes $(X_\zs{i,n})_\zs{n\le \nu}$ and $(X_\zs{i,n})_\zs{n>\nu}$ have the following transition densities in $\bbr^{p}$, $\theta_i\in\Theta_\zs{i}$ 
 \begin{equation}
 \label{dens-f-f^*-1}
 f^{*}_\zs{i}(y\vert x)=
 \frac{\exp\set{-\frac{\vert \eta^{*}_\zs{i}(y,x)\vert^{2}}{2}}}{(2\pi)^{p/2}\sqrt{\det (G_\zs{i}(x))}}
\quad\mbox{and}\quad
 f_\zs{i,\theta_i}(y\vert x)= \frac{\exp\set{-\frac{\vert \eta_\zs{\theta_i,i}(y,x)\vert^{2}}{2}}}{(2\pi)^{p/2}\sqrt{\det (G_\zs{i}(x))}}
\,,
\end{equation}
where $\eta^{*}_\zs{i}(y,x)=G^{-1/2}_\zs{i}(x)(y-\theta^{*}_\zs{i} x)$, 
 $\eta_\zs{\theta_i,i}(y,x)=G^{-1/2}_\zs{i}(x)(y-\theta_i x)$,
 $\theta_i\in\Theta_\zs{i}$ and
 $$
 G_\zs{i}(x)=\EV\, [B_\zs{i,1}xx'B'_\zs{i,1}]+Q^{*}_\zs{i}
=Q_\zs{i}\vect(xx')+Q^{*}_\zs{i}\,.
$$
Therefore, in this case,
$$
g_\zs{i}(\theta_i,y,x)=\log\frac{f_\zs{i,\theta_i}(y\vert x)}{f^{*}_\zs{i}(y\vert x)}
=y'G^{-1}_\zs{i}(x)(\theta_i-\theta^{*}_\zs{i})x
+
\frac{x'(\theta^{*}_\zs{i})^{'}G^{-1}_\zs{i}(x)\,\theta^{*}_\zs{i}x-x'\theta_i'\,G^{-1}_\zs{i}(x)\,\theta_i\,x}{2}
\,.
$$
\noindent
Now  we set
\begin{equation}
\label{L-1-Ex}
\L_\zs{i,n}
=\varrho^{opt}
\sum^{n-1}_\zs{k=0}\,
(1-\varrho^{opt})^{k}\,
 \int_\zs{\Theta_\zs{i}}\, 
e^{\sum^{n}_\zs{j=k+1}g_\zs{i}(\theta_i,X_\zs{i,j},X_\zs{i,j-1})}
  W_\zs{i}(\rm{d}\theta_i)
\end{equation}
and
\begin{equation}
\label{L-1-Ex-22}
\wh{\L}_\zs{i,n}
=
\varrho^{opt}
\sum^{n-1}_\zs{k=0}\,
(1-\varrho^{opt})^{k}
 \sup_\zs{\theta\in\Theta_\zs{i}}\, 
e^{\sum^{n}_\zs{j=k+1}g_\zs{i}(\theta_i,X_\zs{i,j},X_\zs{i,j-1})}
 \,,
\end{equation}
where
$\varrho^{opt}$
is defined in
\eqref{sec:Up-Bnd.1}.
\noindent
 The  random $N\times N$ matrix \eqref{ST-RdMatr}  has the following form
\begin{equation}
\label{ST-RdMatr1}
<\U_\zs{n}>_\zs{i,j}
=\frac{\L_\zs{i,n}}{\wh{\L}_\zs{j,n}}
\quad\mbox{if $i\neq j$ and}\quad
<\U_\zs{n}>_\zs{i,i}=\frac{\L_\zs{i,n}}{\left(1-\varrho^{opt}\right)^{n}}
\,,
\end{equation}
and the corresponding change detection-identification procedure $\delta_\beta^{opt}=(T_\beta^{opt}, d_\beta^{opt})$ is defined by \eqref{seq-prs-opt}-\eqref{Topt} with
the threshold matrix $A=A_\beta^{opt}$ given by \eqref{sec:Cnrsk.7}.

As shown in \cite{PerTar-JMVA2019}, conditions $(\C_\zs{1})$ and $(\C_\zs{2}(r))$ hold
for any $r>0$.  Moreover, one can calculate directly that
 $$
\overline{J}_\zs{i}(\theta_i)=
\frac{1}{2}\,
\EV\,\brcs{\varsigma_\zs{i,\theta_i}'(\theta_i-\theta^{*}_\zs{i})'G^{-1}_\zs{i}(\varsigma_\zs{i,\theta_i})(\theta_i-\theta^{*}_\zs{i})\varsigma_\zs{i,\theta_i}}
$$
and
$$
\overline{J}^{*}_\zs{i}(\theta_i)=
-\frac{1}{2}\,
\EV\, \brcs{
(\varsigma^{*}_\zs{i})'(\theta_i
-\theta^{*}_\zs{i})'G^{-1}_\zs{i}(\varsigma^{*}_\zs{i})(\theta_i-\theta^{*}_\zs{i})\varsigma^{*}_\zs{i}}
\,.
$$
\noindent
To check the condition $(\C_\zs{3})$
denote $\overline{\theta}_\zs{i}=\theta^{*}_\zs{i}-u$. It can be easily shown that for $1\le i\le N$ 
\begin{align*}
g_\zs{i}(u,y,x)
&-\overline{J}^{*}_\zs{i}(u)
=
-
\sum^{p}_\zs{s_\zs{1},s_\zs{2}=1}
<\overline{\theta}_\zs{i}>_\zs{s_\zs{1},s_\zs{2}}
\sum^{p}_\zs{k=1}
D^{(i,1)}_\zs{s_\zs{1},s_\zs{2},k}(x)<\eta^{*}_\zs{i}(y,x)>_\zs{k}\\[2mm]
&
-\frac{1}{2}
\sum^{p}_\zs{s_\zs{1},s_\zs{2},s_\zs{3},s_\zs{4}=1}
<\overline{\theta}_\zs{i}>_\zs{s_\zs{1},s_\zs{2}}
<\overline{\theta}_\zs{i}>_\zs{s_\zs{3},s_\zs{4}}
D^{(i,2)}_\zs{s_\zs{1},s_\zs{2},s_\zs{3},s_\zs{4}}(x)
\,,
\end{align*}
where
 $D^{(i,1)}_\zs{s_\zs{1},s_\zs{2},k}(x)=<G^{-1/2}_\zs{i}(x)>_\zs{s_\zs{1},k}<x>_\zs{s_\zs{2}}$,
 $D^{(i,2)}_\zs{s_\zs{1},s_\zs{2},s_\zs{3},s_\zs{4}}(x)=\b_\zs{s_\zs{1},s_\zs{2},s_\zs{3},s_\zs{4}}(x)-\EV^{*}[\b_\zs{s_\zs{1},s_\zs{2},s_\zs{3},s_\zs{4}}(\varsigma^{*}_\zs{i})]$ and
 $$
 \b_\zs{s_\zs{1},s_\zs{2},s_\zs{3},s_\zs{4}}(x)=<x>_\zs{s_\zs{2}}<x>_\zs{s_\zs{4}}<G^{-1}_\zs{i}(x)>_\zs{s_\zs{1},s_\zs{3}}
 \,.
 $$
\noindent
Note that for any $z\in\bbr^{p}$ and $\vert z\vert =1$
$$
z'G_\zs{i}(x)z=\vert x\vert^{2}(\vect(\wt{x}z'))'Q_\zs{i}\vect(\wt{x}z')+z'Q^{*}_\zs{i}z\,,
$$
where $\wt{x}=x/\vert x\vert$. Taking into account that the matrices $Q_\zs{i}$ are positive definite, we obtain that for some $\c_\zs{*}>0$
\begin{equation}
\label{up-G-11-bnd}
\max_\zs{1\le i\le N}\,
\vert G^{-1}_\zs{i}(x)\vert \le \frac{\c_\zs{*}}{1+\vert x\vert^{2}}
\,.
\end{equation}
\noindent
Therefore,  the functions $D^{(i,1)}_\zs{s_\zs{1},s_\zs{2},k}(x)$ and  $\b_\zs{s_\zs{1},s_\zs{2},s_\zs{3},s_\zs{4}}(x)$ are bounded.
Moreover, as shown in \cite{PerTar-JMVA2019} (Example 1), the Lyapunov function for the process \eqref{sec:Ex.6-0n} (with $\nu=+\infty$)
has the form 
$$
\V(x)=\v^{*}\left( 1+ (x^{\prime}Tx)^{\delta}\right)
$$
for some constant $\v^{*}\ge 1$, a fixed matrix $T$ and any $\delta>0$. Since
in this case
$
\max_\zs{1\le i\le N}\,
\sup_\zs{j\ge  1}\,
\EV^{*}[X^{2}_\zs{i,j}]
<\infty,
$
we obtain that condition $(\C^{*}_\zs{3}(\q))$ holds true for any $\q>0$.
Now, taking into account that under $\Pb^{*}$ the random vectors $(\eta^{*}_\zs{i}(X_\zs{i,j},X_\zs{i,j-1}))_\zs{1\le i\le N\,,\,j\ge 1}$ are i.i.d. $(0,\cI_\zs{p})$ 
Gaussian ($\cI_\zs{p}$ is the unity matrix in $\bbr^{p}$), we 
obtain the condition $(\C_\zs{3}(q))$ for any $q>0$ using Proposition~\ref{Pr.sec:Mrk.1-1}.
 Therefore, Theorems \ref{Th.sec:Cnrsk.2} and  \ref{Th.sec:Cnrsk.2-*}
imply that the sequential procedure $\delta_\beta^{opt}$ defined in \eqref{seq-prs-opt}-\eqref{Topt} 
is asymptotically optimal and robust in the pointwise and minimax senses for any compact sets $\Theta_\zs{i}\subset \Theta^{st}_\zs{i}$ and for any $r>0$. 

\subsection{Example 2:  Change in the correlation coefficients of autoregressive models} \label{ssec:Ex2}

Consider the problem of detecting the change of the correlation coefficient in the $p$th order AR process  which in the $i$th stream satisfies the recursion 
\begin{equation}\label{sec:Ex.7}
X_\zs{i,n} =\a^{(n)}_\zs{i,1}\,X_\zs{i,n-1}+\ldots+\a^{(n)}_\zs{i,p}\,X_\zs{i,n-p}+w_\zs{i,n}, \quad \text{for $n\ge 1$},
\end{equation}
where $\a^{(n)}_\zs{i,l}=\theta^{*}_\zs{l}\Ind{n \le \nu}+\theta_\zs{i,l}\Ind{n > \nu}$
and  $(w_\zs{i,n})_\zs{n\ge 1}$ are i.i.d. Gaussian random variables 
with $\EV\,[w_\zs{i,1}]=0$, $\EV\,[w^{2}_\zs{i,1}]=1$.  
In the sequel, we use the notation  $\theta^{*}_\zs{i}=(\theta^{*}_\zs{i,1},\ldots,\theta^{*}_\zs{i,p})^\prime$ and 
$\theta_\zs{i}=(\theta_\zs{i,1},\ldots,\theta_\zs{i,p})^\prime$. Hereafter $\prime$ denotes the transposition operation.
The corresponding conditional densities $X_\zs{i,n}\vert X_\zs{i,n-1},\ldots,X_\zs{i,n-p}$ for  $\nu\le n$ and $\nu>n$, $\theta_i\in\Theta_\zs{i}$, are
 \begin{equation}
 \label{dens-AR-p-1}
 f^{*}_\zs{i}(y\vert x)=
 \frac{1}{(2\pi)^{p/2}}
 \,\exp\set{-\dfrac{
(\eta^{*}_\zs{i}(y,x))^{2}}{2}}
\quad\mbox{and}\quad
 f_\zs{i,\theta_i}(y\vert x)=
  \frac{1}{(2\pi)^{p/2}}
 \,\exp\set{-\dfrac{
(\eta_\zs{\theta_i,i}(y,x))^{2}}{2}}
\,,
\end{equation}
where $\eta^{*}_\zs{i}(y,x)=y-(\theta^{*}_\zs{i})^{\prime} x$ and $\eta_\zs{\theta_i,i}(y,x)=y-(\theta_\zs{i})^{\prime} x$.
\noindent
Therefore,  for any $\theta_i\in\bbr^{p}$, $y\in\bbr$ and  
$x=(x_\zs{1},\ldots,x_\zs{p})^{\prime}\in\bbr^{p}$
\begin{equation}\label{func-g-Coef-nn-1}
g_\zs{i}(\theta_i,y,x)=
\log\frac{f_\zs{i,\theta_i}(y\vert x)}{f^{*}_\zs{i}(y\vert x)}
=
y_\zs{1}(\theta_i-\theta^{*}_\zs{i})^{\prime}x
+\frac{((\theta^{*}_\zs{i})^{\prime} x)^{2}-(\theta_i^{\prime}x)^{2}}{2}
\,.
\end{equation}
 The process \eqref{sec:Ex.7} is not Markov, but the $p$-dimensional processes
\begin{equation}\label{sec:Ex.7-1n}
\Phi_\zs{i,n}=(X_\zs{i,n},\ldots,X_\zs{i,n-p+1})^{\prime}\in\bbr^{p}, \quad i=1,\dots,N
\end{equation}
are Markov. Now, for any $\vartheta=(\vartheta_\zs{1},\ldots,\vartheta_\zs{p})\in\bbr^{p}$ we difine
\[
\Lambda(\vartheta)= \begin{pmatrix}
\vartheta_\zs{1} & \vartheta_2 & \dots & \vartheta_\zs{p}\\
1 & 0 &  \dots & 0\\
\vdots & \vdots &  \ddots & \vdots \\
0 & 0 &\dots 1&0 
\end{pmatrix} \, .
\]
Using this matrix it is easy to show that the  processes $(\Phi_\zs{i,n})_\zs{n\le \nu}$ and $(\Phi_\zs{i,n})_\zs{n>\nu+p}$
satisfy the following stochastic linear equations:
\begin{equation}\label{sec:Ex.8}
\Phi_\zs{i,n}=\Lambda_\zs{*,i} \Phi_\zs{i,n-1}+\wt{w}_\zs{i,n} \quad \mbox{for $n\le \nu$} \quad \text{and} \quad
\Phi_\zs{i,n}=\Lambda_\zs{i}\Phi_\zs{i,n-1}+\wt{w}_\zs{i,n} \quad\mbox{for $n>\nu$} \, ,
\end{equation}
where $\Lambda_\zs{*,i}=\Lambda(\theta^{*}_\zs{i})$, $\Lambda_\zs{i}=\Lambda(\theta_\zs{i})$ and $\wt{w}_\zs{i,n}=(w_\zs{i,n},0,\ldots,0)'\in\bbr^{p}$.
Obviously,
\[
\EV[\wt{w}_\zs{n}\,\wt{w}^{\prime}_\zs{n}] =B=
\begin{pmatrix}
1 & \dots & 0\\
\vdots & \ddots & \vdots\\
0 & \dots & 0
\end{pmatrix} .
\]
\noindent
Assume that 
all eigenvalues of the matrices $\Lambda_\zs{*,i}=\Lambda_\zs{i}(\theta^{*}_\zs{i})$ in modules are less than $1$ and that $\theta_\zs{i}$ belongs to the set 
\begin{equation}
\label{set-Theta-i}
\Theta^{st}_\zs{i}=\{\vartheta\in \bbr^{p}\,:\, \max_\zs{1\le j\le p}\,\vert
\e_\zs{j}(\Lambda_\zs{i}(\vartheta))\vert<1\}
\setminus \,
\{\theta^{*}_\zs{i}\}
\,,
\end{equation}
where $\e_\zs{j}(\Lambda)$ is the $i$th eigenvalue of the matrix $\Lambda$.
In this case, the processes \eqref{sec:Ex.8} have the ergodic distributions defined by the random vectors 
$$ 
\varsigma^{*}_\zs{i}=\sum_\zs{l\ge 1}\, (\Lambda_\zs{*,i})^{l-1} \wt{w}_\zs{i,l} 
\quad\mbox{and}\quad
\varsigma_\zs{i,\theta}=\sum_\zs{l\ge 1} (\Lambda_\zs{i})^{l-1} \wt{w}_\zs{i,l}
$$
which are  $(0,\F^{*})$ and $(0,\F_\zs{i}(\vartheta))$ Gaussian vectors in $\bbr^{p}$, where 
$$
\F^{*}
=\sum_\zs{n\ge 0} (\Lambda_\zs{*,i})^{n}\, B\, (\Lambda'_\zs{*,i})^{n}
\quad\mbox{and}\quad
\F_\zs{i}(\vartheta)=\sum_\zs{n\ge 0} \Lambda^{n}_\zs{i} B (\Lambda'_\zs{i})^{n}
\,.
$$
\noindent
Now  we set
\begin{equation}
\label{L-1-Ex-EX-2}
\L_\zs{i,n}
=\varrho^{opt}
\sum^{n-1}_\zs{k=0}\,
(1-\varrho^{opt})^{k}\,
 \int_\zs{\Theta_\zs{i}}\, 
e^{\sum^{n}_\zs{j=k+1}g_\zs{i}(\vartheta,X_\zs{i,j},\Phi_\zs{i,j-1})}
  W_\zs{i}(\rm{d}\vartheta)
\end{equation}
and
\begin{equation}
\label{L-1-Ex-22-EX2-2}
\wh{\L}_\zs{i,n}
=
\varrho^{opt}
\sum^{n-1}_\zs{k=0}\,
(1-\varrho^{opt})^{k}
 \sup_\zs{\theta\in\Theta_\zs{i}}\, 
e^{\sum^{n}_\zs{j=k+1}g_\zs{i}(\vartheta,X_\zs{i,j},\Phi_\zs{i,j-1})}
 \,,
\end{equation}
where
$\varrho^{opt}$
is defined in
\eqref{sec:Up-Bnd.1}.
\noindent
The  random $N\times N$ matrix \eqref{ST-RdMatr} 
has the following form
\begin{equation}
\label{ST-RdMatr-EX2}
<\U_\zs{n}>_\zs{i,j}
=\frac{\L_\zs{i,n}}{\wh{\L}_\zs{j,n}}
\quad\mbox{if $i\neq j$}\quad \text{and} \quad
<\U_\zs{n}>_\zs{i,i}=\frac{\L_\zs{i,n}}{\left(1-\varrho^{opt}\right)^{n}}
\,.
\end{equation}
and the corresponding change detection-identification procedure $\delta_\beta^{opt}=(T_\beta^{opt}, d_\beta^{opt})$ is defined in \eqref{seq-prs-opt}-\eqref{Topt} with
the threshold matrix $A=A_\beta^{opt}$ given by \eqref{sec:Cnrsk.7}.

As shown in \cite{PerTar-JMVA2019}, conditions $(\C^{'}_\zs{1})$ and $(\C_\zs{2}(r))$ hold
for any $r>0$ and any compact sets $ \Theta_\zs{i}\subset \Theta^{st}_\zs{i}$  for the  function $I_\zs{i,j}$ defined in \eqref{KL-Inf-1} with
\begin{equation}
\label{KLB-IBf--12sc}
\overline{J}_\zs{i}(\vartheta)=
\frac{1}{2}\,
(\vartheta-\theta^{*}_\zs{i})'
\F(\vartheta)
(\vartheta-\theta^{*}_\zs{i})
\quad\mbox{and}\quad
\overline{J}_\zs{i}(\vartheta)=
-
\frac{1}{2}\,
(\vartheta-\theta^{*}_\zs{i})'
\F^{*}
(\vartheta-\theta^{*}_\zs{i})\,.
\end{equation}
It should be noted that in the scalar case, i.e., when $p_\zs{1}=\ldots= p_\zs{N}=1$, 
$$
\overline{J}_\zs{i}(\vartheta)=
\frac{(\vartheta-\theta^{*}_\zs{i})^{2}}{2(1-\vartheta^{2})}
\quad\mbox{and}\quad
\overline{J}_\zs{i}(\vartheta)=
-
\frac{(\vartheta-\theta^{*}_\zs{i})^{2}}{2(1-(\vartheta^{*}_\zs{i})^{2})}\,.
$$
\noindent
Write $\overline{\theta}_\zs{i}=\theta^{*}_\zs{i}-u$. To check the condition $(\C_\zs{3})$ direct calculations show that for $1\le i\le N$ 
$$
g_\zs{i}(u,y,x)
-\overline{J}^{*}_\zs{i}(u)
=
-\overline{\theta}^{\prime}_\zs{i}\eta_\zs{i}(y,x)
-
\frac{1}{2}
\overline{\theta}^{\prime}_\zs{i}
\left(
xx^{\prime}
-
\EV [\varsigma^{*}_\zs{i,1} (\varsigma^{*}_\zs{i,1})^{\prime}] 
\right)
\overline{\theta}_\zs{i}
\,.
$$
\noindent
 Therefore, taking into account that
in this case 
$
\max_\zs{1\le i\le N}\,
\sup_\zs{j\ge  1}\,
\EV^{*}\vert X_\zs{i,j}\vert^{\q}
<\infty
$
for any $\q>0$ we obtain that condition $(\C^{*}_\zs{3}(\q))$ holds for any $\q>0$.
Now, taking into account that under the probability measure $\Pb^{*}$ the random variables 
$(\eta^{*}_\zs{i}(X_\zs{i,j},\Phi_\zs{i,j-1}))_\zs{1\le i\le N\,,\,j\ge 1}$ are i.i.d. $\Nc(0,1)$ 
we obtain the condition $(\C_\zs{3}(q))$  for any $q>0$ using  Proposition~\ref{Pr.sec:Mrk.1-1}.
Therefore, Theorems \ref{Th.sec:Cnrsk.2} and  \ref{Th.sec:Cnrsk.2-*}
imply that the sequential procedure $\delta_\beta^{opt}$ is asymptotically optimal and robust in the pointwise and minimax senses for any compact sets $\Theta_\zs{i}\subset \Theta^{st}_\zs{i}$
and for any $r>0$.

\section{Application to epidemics detection and localization} \label{sec:Epid} 

\subsection{Near optimality}\label{ssec:AO}

We begin with considering the epidemiological statistical models proposed in  \cite{BaronChoudharyYu2013}.
Assume that for any $1\le i\le N$ the observations $(X_\zs{i,n})_\zs{1\le n\le \nu}$ and 
 $(X_\zs{i,n})_\zs{n> \nu}$ are homogenous  Markov processes with the values in the finite space 
 $(\Xc,\mu)$, $\Xc_\zs{i}=\{0,\ldots,D\}$ and $\mu{\{0\}}=\ldots=\mu{\{D\}}=1$. 
In this model, the conditional
 $X_\zs{i,n}\vert X_\zs{i,n-1}$
  densities  for $n\le \nu$
and for $n>\nu$ are
defined respectively as
\begin{equation}
\label{ep-dn-1}
f^{*}_\zs{i}(y \vert x)=
\dbinom{x}{y}
(p^{*}_\zs{i})^{x-y} (1-p^{*}_\zs{i})^{y}
\quad\mbox{and}\quad
f_\zs{i,\theta_i}(y \vert x)=
\dbinom{x}{y}
\theta_i^{x-y} (1-\theta_i)^{y}
\,,
\end{equation}
where $0<p^{*}_\zs{i}<1$ and $\theta_i\in\Theta_\zs{i}\subset [0,1]$. The probabilities $p^{*}_\zs{i}$ are non-epidemic (normal)  infection rates
and the $\Theta_\zs{i}$ are the sets of epidemic values of the infection parameters $\theta_i$.  In this case, the functions $g_\zs{i}$ defined in
\eqref{sec:Mrk.5}  for any $0<\theta_i<1$, $x,y\in\Xc$ have the following forms
\begin{equation}
\label{fun-g-ep-1}
g_\zs{i}(\theta_i,y,x)=\log \frac{f_\zs{i,\theta_i}(y|x)}{f^{*}_\zs{i}(y|x)}
=(x-y)\log\frac{\theta_i}{p^{*}_\zs{i}}+y\log\frac{1-\theta_i}{1-p^{*}_\zs{i}}, \quad i=1,\dots,N.
\,.
\end{equation}
So the functions
\eqref{sec:Mrk.6} are
\begin{equation*}
J_\zs{i}(\theta_i,x)
=
\sum^{x}_\zs{y=0}\,
\left(
(x-y)\log\frac{\theta_i}{p^{*}_\zs{i}}
+y\log\frac{1-\theta_i}{1-p^{*}_\zs{i}}
 \right)
 \dbinom{x}{y}
\theta_i^{x-y} (1-\theta_i)^{y}
\end{equation*}
and
\begin{equation*}
J^{*}_\zs{i}(\theta_i,x)
=
\sum^{x}_\zs{y=0}\,
\left(
(x-y)\log\frac{\theta_i}{p^{*}_\zs{i}}
+y\log\frac{1-\theta_i}{1-p^{*}_\zs{i}}
 \right)
 \dbinom{x}{y}
\, (p^{*}_\zs{i})^{x-y} (1-p^{*}_\zs{i})^{y}
 \,.
\end{equation*}

One can check directly that the set $C=\{0\}$ is an accesible atom for the Markov chains  $(X_\zs{i,n})_\zs{n> \nu}$.  Obviously, if $X_\zs{i,n-1}=0$ then $X_\zs{i,n}=0$ almost surely.
Define the Markov time 
 $$
 \tau=\inf\{n\ge 1\,:\, X_\zs{i,n}=0\}
 \,.
 $$
 If $X_\zs{i,0}=0$ then $\tau=1$, i.e.,  $\EV^{\theta_i}_\zs{i,x=0}\,[\tau]=1$.
 Therefore, for any $0<\theta_i<1$ the chain is ergodic with the ergodic distribution $\lambda^{\theta_i}(\Gamma)=\Chi_\zs{\{0\in\Gamma\}}$
 for any $\Gamma\subseteq\Xc$ (point measure). See, e.g., Theorems 10.2.1 and 10.2.2 in \cite{MeTw93}. In this case, for any compact sets $\Theta_\zs{i}\subset\Theta$
 $$
 \min_\zs{1\le i\le N}\,
 \inf_\zs{\theta_i\in\Theta_\zs{i}}\,
 \inf_\zs{x,y\in C}\,f_\zs{i,\theta_i}(y|x)
 =1\,.
$$
Let us select
\begin{equation}
\label{func-V-1}
\V(x)=\V_\zs{*}e^{\gamma x}
\,,
\end{equation}
where $\gamma>0$ and $\V_\zs{*}\ge 1$.
For any $x\ge 1$ we have
\begin{align*}
\frac{\EV^{\theta_i}_\zs{x}[\V(X_\zs{1})]}{\V(x)}
&=\sum^{x}_\zs{y=0}\,
\frac{\V(y)}{\V(x)}\,
\dbinom{x}{y}
\theta_i^{x-y} (1-\theta_i)^{y}\\[3mm]
&=(1-\theta_i)^{x}
+
\sum^{x-1}_\zs{y=0}\,
\frac{\V(y)}{\V(x)}\,
\dbinom{x}{y}
\theta_i^{x-y} (1-\theta_i)^{y}\\[3mm]
&\le 
(1-\theta_i)^{x}
+e^{-\gamma}
\sum^{x-1}_\zs{y=0}\,
\dbinom{x}{y}
\theta_i^{x-y} (1-\theta_i)^{y}
=(1-e^{-\gamma})(1-\theta_i)^{x}
+e^{-\gamma}
\\[3mm]
&\le 
(1-e^{-\gamma})(1-\theta_i)
+e^{-\gamma}
\,.
\end{align*}
So, if we take $\gamma=\log 2$ in \eqref{func-V-1}, we obtain that for $1\le i\le N$ and $\theta_i\in\Theta_\zs{i}$
$$
\frac{\EV^{\theta_i}_\zs{x}[\V(X_\zs{i,1})]}{\V(x)}\le \frac{2-\theta_i}{2}
\le
1-\rho
\quad\mbox{with}\quad
\rho=\min_\zs{1\le i\le N}\,
\inf_\zs{\theta_i\in\Theta_\zs{i}}\,
\frac{\theta_i}{2}
\,.
$$

\noindent
By Theorem~ A1 in \cite{PerTar-JMVA2019}, the Markov chain $(X_\zs{i,n})_\zs{n> \nu}$
is uniformly geometric ergodic and for some positive constants $\kappa^{*}$ and $R^{*}$
$$
\sup_\zs{n\ge 0}\,
e^{\kappa^{*} n}\,
\max_\zs{1\le i\le N}\,
\sup_\zs{x\in\Xc_\zs{i}}\,
\sup_\zs{\theta_i\in\Theta_\zs{i}}
\sup_\zs{0\le g\le \V}
\frac{1}{\V(x)}
|
\EV^{\theta}_\zs{i,x}\,\,[g(X_\zs{i,n})]
-
g(0)
|
\le R^{*}
\,.
$$
Therefore, Theorem~\ref{Th.sec:Mrk.1} implies conditions $(\A_\zs{1}$) and $(\A_\zs{2}(r))$ with
$\overline{J}_\zs{i}(\theta_\zs{i})=J_\zs{i}(\theta_i,0)=0$ and $\overline{J}^{*}_\zs{j}(\theta_\zs{j})=J^{*}_\zs{i}(\theta_i,0)=0$ for all $r>0$.
This means that we cannot use the procedures \eqref{WSR-def} for this problem directly. However, 
in practice the values of the observations $X_\zs{i,n}$ are sufficiently large, i.e., $D\to \infty$, and usually
the number of the infected populations is not too large, i.e.,
 $X_\zs{i,n}\ge \epsilon D$ for some $0<\epsilon<1$. So it is more natural to modify the initial model and study the limiting model when $D$ is sufficiently large. 
 Note that in this case observations in the binomial models \eqref{ep-dn-1} can be represented as  
 $$
 X_\zs{i,n}=(1-\vartheta)  X_\zs{i,n-1}+\sum^{X_\zs{i,n-1}}_\zs{j=1}\,(\eta_\zs{n,j}-1+\vartheta)
 \,,
 $$
 where $(\eta_\zs{n,j})_\zs{j\ge 1}$ is i.i.d. sequence of Bernoulli random variables with 
 $\Pb(\eta_\zs{n,j}=1)=1-\vartheta$ and
 independent from $X_\zs{i,n-1}$ and where $\vartheta=\theta_i$ and $p_i^*$ in the post-change and pre-change modes, respectively.  Using  the Gaussian approximation for the last sum
 $$
\frac{1}{\sqrt{X_\zs{i,n-1}}}\,
 \sum^{X_\zs{i,n-1}}_\zs{j=1}\,(\eta_\zs{n,j}-1+\vartheta)
 \sim
 \Nc(0, \sigma^{2}_{\vartheta})\,,\quad
 \sigma^{2}_\zs{\vartheta}=\vartheta(1-\vartheta)
 \,,
 $$
 we obtain the following model
 $$
 X_\zs{i,n}=(1-\vartheta)  X_\zs{i,n-1}+\sigma_\zs{\vartheta}\sqrt{|X_\zs{i,n-1}|}\xi_\zs{i,n}
 \,,
 $$
 where $(\xi_\zs{i,n})_\zs{n\ge 1}$ is the sequence of i.i.d. normal $\Nc(0,1)$ random variables. 
 
 Thus, in place of the original Bernoulli model we will use the following model: 
 the observations $X_\zs{i,k}$ before change are defined as 
 \begin{equation}
 \label{ep-15-05-1}
 X_\zs{i,n}=(1-p^{*}_\zs{i})  X_\zs{i,n-1}+\sigma^{*}_\zs{i}\sqrt{\vert X_\zs{i,n-1}\vert}\,\xi_\zs{i,n}\,,\quad
 \sigma^{*}_\zs{i}=\sqrt{p^{*}(1-p^{*})}\,,
 \end{equation}
and after change as
\begin{equation}
 \label{ep-15-05-2}
 X_\zs{i,n}=(1-\theta_i)  X_\zs{i,n-1}+\sigma_\zs{\theta_i}\sqrt{\vert X_\zs{i,n-1}\vert}\,\xi_\zs{i,n}\,,
 \end{equation}
where $\theta_i\in\Theta_\zs{i}\subset [0,1]$ and $(\xi_\zs{i,n})_\zs{n\ge 1}$ are i.i.d. $\Nc(0,1)$ random variables. 
In this case, the spaces $(\Xc_\zs{i},\Bc_\zs{i},\mu_\zs{i})$ are:  $\Xc=\bbr_\zs{*}=\bbr\setminus \{0\}$, $\Bc_\zs{i}=\Bc(\bbr_\zs{*})$
is the Borel field
 and 
$\mu_\zs{i}=\mu$ is the Lebesgue measure on $\Bc(\bbr_\zs{*})$.
Obviously,
\begin{equation}
\label{dntes-1}
f^{*}_\zs{i}(y\vert x)=\frac{1}{\sigma^{*}_\zs{i}\sqrt{2\pi\vert x\vert}}\,\exp\set{-\frac{[y-(1-p^{*})x]^{2}}{2(\sigma^{*}_\zs{i})^{2}\vert x\vert}}
\quad\mbox{and}\quad
f_\zs{\theta_i}(y\vert x)=\frac{1}{\sigma_\zs{\theta_i}\sqrt{2\pi\vert x\vert}}\,\exp\set{-\frac{[y-(1-\theta_i)x]^{2}}{2\sigma^{2}_\zs{\theta_i}\vert x\vert}} \, .
\end{equation}
Using definitions \eqref{sec:Mrk.6-h-def} and \eqref{sec:Mrk.6} we obtain that for $x\in\bbr_\zs{*}$
\begin{equation}
\label{g-1-16-05-2-0}
g_\zs{i}(\theta_i,y,x)
=\log\frac{\sigma^{*}_\zs{i}}{\sigma_\zs{\theta_i}}
+\frac{(\eta^{*}_\zs{i}(y,x))^{2}}{2}
-
\frac{(\eta_\zs{\theta_i,i}(y,x))^{2}}{2}
\,,
\end{equation}
where 
$$
\eta^{*}_\zs{i}(y,x)=\frac{y-(1-p^{*})x}{\sigma^{*}_\zs{i}\sqrt{\vert x\vert}}
\quad\mbox{and}\quad
\eta_\zs{\theta_i,i}(y,x)=
\frac{y-(1-\theta_i)x}{\sigma^{*}_\zs{\theta_i}\sqrt{\vert x\vert}} .
$$
In this case,
\begin{equation}
\label{L-1-Ex-EP-0}
\L_\zs{i,n}
=\varrho^{opt}
\sum^{n-1}_\zs{k=0}\,
(1-\varrho^{opt})^{k}\,
 \int_\zs{\Theta_\zs{i}}\, 
e^{\sum^{n}_\zs{j=k+1}g_\zs{i}(\theta_i,X_\zs{i,j},X_\zs{i,j-1})}
  W_\zs{i}(\rm{d}\theta_i)
\end{equation}
and
\begin{equation}
\label{L-1-Ex-22-EP-1}
\wh{\L}_\zs{i,n}
=
\varrho^{opt}
\sum^{n-1}_\zs{k=0}\,
(1-\varrho^{opt})^{k}
 \sup_\zs{\theta\in\Theta_\zs{i}}\, 
e^{\sum^{n}_\zs{j=k+1} g_\zs{i}(\theta_i,X_\zs{i,j},X_\zs{i,j-1})}
 \,,
\end{equation}
where $\varrho^{opt}$ is defined in \eqref{sec:Up-Bnd.1}. The  elements of the random $N\times N$ matrix \eqref{ST-RdMatr} have the following form
\begin{equation}
\label{ST-RdMatr-EP-33}
<\U_\zs{n}>_\zs{i,j}
=\frac{\L_\zs{i,n}}{\wh{\L}_\zs{j,n}}
\quad\mbox{if $i\neq j$}\quad \text{and} \quad 
<\U_\zs{n}>_\zs{i,i}=\frac{\L_\zs{i,n}}{\left(1-\varrho^{opt}\right)^{n}}
\, ,
\end{equation}
and the corresponding change detection-identification procedure $\delta_\beta^{opt}=(T_\beta^{opt}, d_\beta^{opt})$ is defined in \eqref{seq-prs-opt}-\eqref{Topt} with
the threshold matrix $A=A_\beta^{opt}$ given by \eqref{sec:Cnrsk.7}.

Let us check conditions $(\C_\zs{1})-(\C_\zs{3})$. To this end, first note that
\begin{equation}
\label{g-1-17-05-tl-g-1}
J_\zs{i}(\theta_i,x)
=\frac{1}{2}
\left(
\log\frac{p^{*}_\zs{i}(1-p^{*}_\zs{i})}{\theta_i(1-\theta_i)}
-1
+\frac{\theta_i(1-\theta_i)}{p^{*}_\zs{i}(1-p^{*}_\zs{i})}
+\frac{(\theta_i-p^{*}_\zs{i})^{2}}{p^{*}_\zs{i}(1-p^{*}_\zs{i})}\vert x\vert
\right)
\end{equation}
and
\begin{equation}
\label{g-J-overline-1}
J^{*}_\zs{i}(\theta_i,x)
=\frac{1}{2}
\left(
\log\frac{p^{*}_\zs{i}(1-p^{*}_\zs{i})}{\theta_i(1-\theta_i)}
+1
-\frac{p^{*}_\zs{i}(1-p^{*}_\zs{i})}{\theta_i(1-\theta_i)}
-\frac{(\theta_i-p^{*}_\zs{i})^{2}}{{\theta_i(1-\theta_i)}}\vert x\vert
\right)
\,.
\end{equation}
\noindent 
Taking into account that the function $z-\log z-1>0$ for all $z>0$ and $z\neq 1$ we obtain that
$\inf_\zs{x} J_\zs{i}(\theta_i,x)>0$ and $\sup_\zs{x} J_\zs{i}(\theta_i,x)<0$  for all $\theta_i\neq p^{*}_\zs{i}$. 
Recall that $p_i^{*}$ and $\theta_i$ are the infection rates, where $p^{*}$ is normal non-epidemic value and $\theta_i$ is epidemic value. So, if 
$\check{p}_\zs{i}$ is the epidemic threshold for the $i$th stream, then 
$\Theta_\zs{i}\subset (\check{p}_\zs{i},1)$
and for
all $p^{*}_\zs{i}<\check{p}_\zs{i}$ 
\begin{equation}
\label{lower-bnd}
\min_\zs{1\le i\le N}\,
\inf_\zs{\theta_i\in\Theta}\,
\inf_\zs{x}
J_\zs{i}(\theta_i,x)>0
\,.
\end{equation}
\noindent
Now we need to check the conditions of 
Theorem \ref{Th.sec:Mrk.1}. First, note that the definition
\eqref{dntes-1}
yields that for any $\B>0$
$$
\min_\zs{1\le i\le N}
\,
\inf_\zs{\theta_i\in\Theta_\zs{i}}
\,
\inf_\zs{\vert y\vert\le \B}
\,
\inf_\zs{\vert x\vert\le \B}
\,
f_\zs{\theta_i}(y\vert x)
>0\,.
$$
From
\eqref{g-1-17-05-tl-g-1}
it is easy to deduce that for any compact sets $\Theta_\zs{i}\subset (\check{p}_\zs{i},1)$
\begin{equation}
\label{Up-bnd-g^*-1}
\g^{*}=
\max_\zs{1\le i\le N}\,
\max_\zs{x\in\bbr}\,
\sup_\zs{\theta_i\in\Theta_\zs{i}}
\frac{J_\zs{i}(\theta_i,x)}{1+\vert x\vert}<\infty\,.
\end{equation}
\noindent
Moreover, note that
$$
\frac{\partial g_\zs{i}(\vartheta,y,x)}{\partial \vartheta}
=\frac{2\vartheta-1}{\vartheta(1-\vartheta)}
-\frac{(y-(1-\vartheta)x)x}{\vartheta(1-\vartheta)\vert x\vert}
+
\frac{(y-(1-\vartheta)x)^{2}}{2\vartheta^{2}(1-\vartheta)^{2}\vert x\vert}
$$
and, therefore, using the definition \eqref{sec:Mrk.6-h-def} we obtain 
\begin{equation}
\label{h-i-x-y-1ubnd}
h(x,y)=\max_\zs{1\le i\le N}
\max_\zs{u\in\Theta_\zs{i}}\,
\vert \partial g_\zs{i}(u,y,x)/
\partial u
\vert
\le 
\check{h}
\left(
1+\vert x\vert
+\vert y\vert
+\frac{y^{2}}{\vert x\vert}
\right)\,,
\end{equation}
where
$
\check{h}=
\max_\zs{1\le i\le N}
\max_\zs{\vartheta\in\Theta_\zs{i}}\,\vartheta^{-2}(1-\vartheta)^{-2}$.
\noindent
Therefore, taking into account, that in this case $f_\zs{i,\vartheta}(y|x)=f_\zs{\vartheta}(y|x)$ for all $1\le i\le N$, we get
for any
$\vartheta\in\Theta_\zs{i}$
\begin{align*}
\overline{h}_\zs{i}(\vartheta,y)&=\overline{h}(\vartheta,y)=
\int_\zs{\bbr}\,h(y,x)\,f_\zs{\vartheta}(y|x)\,\d y
\le 
\check{h}
\left(
1+\vert x\vert
+
\int_\zs{\bbr}\,\vert y\vert\,f_\zs{\vartheta}(y|x)\,\d y
+
\frac{1}{\vert x\vert}
\int_\zs{\bbr}\,y^{2}\,f_\zs{\vartheta}(y|x)\,\d y
\right)\\[3mm]
&\le 
\check{h}
\left(
1+\vert x\vert
+
(1-\vartheta)
\vert x\vert
+\sigma_\zs{\vartheta}\sqrt{\vert x\vert}
+
(1-\vartheta)^{2}
\vert x\vert
+
\sigma^{2}_\zs{\vartheta} \vert x\vert
\right)
\le 4
\check{h}
(
1+\vert x\vert)
\,.
\end{align*}

To check the conditions $({\rm C}1.2)$ we set
$$
\V(x)=\V^{*}(1+\vert x\vert)
\quad\mbox{and}\quad
\V^{*}=\g^{*}+4\check{h}
\,.
$$
For the model \eqref{ep-15-05-2} we have
$$
\EV^{\theta_i}_\zs{i,x}\,[\V(X_\zs{i,1})]\le 
\V^{*}
\left(
1
+
(1-\theta_i)
\vert x\vert
+\sigma_\zs{\theta_i}\sqrt{\vert x\vert}\EV\vert\xi_\zs{1}\vert
\right)
\,,
$$
i.e.,
$$
\limsup_\zs{\B\to\infty}\,
\sup_\zs{\vert x\vert>\B}\,
\max_\zs{1\le i\le N}\,
\sup_\zs{\theta_i\in\Theta_\zs{i}}\,
\frac{\EV^{\theta_i}_\zs{i,x}\,[\V(X_\zs{i,1})]}{\V(x)}
\le 1-
\min_\zs{1\le i\le N}\,
\inf_\zs{\theta_i\in\Theta_\zs{i}}
\theta_i
\,.
$$
Therefore, there exists $\B>0$ such that for all $\vert x\vert>\B$ 
$$
\max_\zs{1\le i\le N}\,
\sup_\zs{\theta_i\in\Theta_\zs{i}}
\,
\EV^{\theta_i}_\zs{i,x}\,[\V(X_\zs{i,1})]\le (1-\rho)\V(x)
\quad\mbox{and}\quad
\rho=
\frac{1}{2}\,
\min_\zs{1\le i\le N}\,
\inf_\zs{\theta_i\in\Theta_\zs{i}}
\theta_i
\,.
$$
\noindent 
Obviously, this inequality  implies condition $(\C_\zs{1})$ with $C=\{x\in\bbr\,:\, \vert x\vert \le \B\}$.
Using Theorem 15.01 in \cite{MeTw93} and  Proposition~\ref{Pr.sec:App.1-1} in the Appendix, it is easy to deduce that the processes
\eqref{ep-15-05-1} and \eqref{ep-15-05-2} are stationary with the ergodic distributions defined by the random variables
$\varsigma^{*}_\zs{i}$ and $\varsigma_\zs{\theta_i,i}$ such that  for any $\q>0$
\begin{equation}
\label{up-bnd-erg-dis-11}
\E\vert \varsigma^{*}_\zs{i}\vert^{\q}
<\infty
\quad\mbox{and}\quad
\E\vert\varsigma_\zs{\theta_i,i}\vert^{\q}
<\infty\,.
\end{equation}
Hence, from \eqref{g-1-17-05-tl-g-1} and \eqref{g-J-overline-1} we get
\begin{equation}
\label{Kul-Back}
\overline{J}_\zs{i}(\theta_i)
=\frac{1}{2}
\left(
\log\frac{p^{*}_\zs{i}(1-p^{*}_\zs{i})}{\theta_i(1-\theta_i)}
-1
+\frac{\theta_i(1-\theta_i)}{p^{*}_\zs{i}(1-p^{*}_\zs{i})}
+\frac{(\theta_i-p^{*}_\zs{i})^{2}}{p^{*}_\zs{i}(1-p^{*}_\zs{i})}\,\EV\vert\varsigma_\zs{\theta_i,i}\vert
\right)
\end{equation}
and 
\begin{equation}
\label{Kul-Back-*} 
\overline{J}^{*}_\zs{i}(\theta_i)
=\frac{1}{2}
\left(
\log\frac{p^{*}_\zs{i}(1-p^{*}_\zs{i})}{\theta_i(1-\theta_i)}
+1
-\frac{p^{*}_\zs{i}(1-p^{*}_\zs{i})}{\theta_i(1-\theta_i)}
-\frac{(\theta_i-p^{*}_\zs{i})^{2}}{{\theta_i(1-\theta_i)}}
\EV\vert \varsigma^{*}_\zs{i}\vert
\right)
\,.
\end{equation}

As far as the condition $(\C_\zs{2})$ is concerned it follows from \eqref{g-1-16-05-2-0} and
\eqref{h-i-x-y-1ubnd} that there exists a constant $\c_\zs{*}>0$ such that for all $x,y\in\bbr_\zs{*}$
$$
\max_\zs{1\le i\le N}\,
\sup_\zs{\theta_i\in\Theta_\zs{i}}
\left(
\vert g_\zs{i}(\theta_i,y,x)\vert+
h_\zs{i}(x,y)
\right)
\le 
\c_\zs{*}
\left(
1+\vert x\vert
+\vert y\vert
+\frac{y^{2}}{\vert x\vert}
\right)\,.
$$
Note that for any $m\ge 1$ there exists a constant  $C_\zs{m}>0$ for which
for any $1\le i\le N$ and $\theta_i\in\Theta_\zs{i}$
$$
\EV_\zs{i,\theta_i}\,\left[ X^{2m}_\zs{i,n}\vert X_\zs{i,n-1}\right] \le C_\zs{m}
(X^{2m}_\zs{i,n-1}+\vert X_\zs{i,n-1}\vert^{m})
\,.
$$
Hence, for any  $m\ge 1$  there exists a constant  $C_\zs{m}>0$ such that
$$
\EV^{\theta_i}_\zs{i,x}\,[g^{2m}_\zs{i}(\theta_i, X_\zs{i,n}, X_\zs{i,n-1})]
\le C_\zs{m}
\left(
1+
\EV^{\theta_i}_\zs{i,x}\,[X^{2m}_\zs{i,n-1}]
+
\EV^{\theta_i}_\zs{i,x}\,[\vert X_\zs{i,n-1}\vert^{m}]
\right)
\le 
C_\zs{m}
\left(
1+
\EV^{\theta_i}_\zs{i,x}\,[X^{2m}_\zs{i,n-1}]
\right)
$$
and
$$
\EV^{\theta_i}_\zs{i,x}\,[h^{2m}_\zs{i}( X_\zs{i,n}, X_\zs{i,n-1})]
\le C_\zs{m}
\left(
1+
\EV^{\theta_i}_\zs{i,x}\,[X^{2m}_\zs{i,n-1}]
\right)
\,.
$$
\noindent
To check the condition $(\C_\zs{3})$ we can obtain directly from
\eqref{Kul-Back-*} that for $1\le i\le N$ 
\begin{align*}
g_\zs{i}(u,y,x)
&-\overline{J}^{*}_\zs{i}(u)
=\frac{1}{2}\left(1
-\frac{p^{*}_\zs{i}(1-p^{*}_\zs{i})}{u(1-u)}\right)((\eta^{*}_\zs{i}(y,x))^{2}-1)
-\frac{(u-p^{*}_\zs{i})^{2}}{2u(1-u)}
\left(\vert x\vert- \E[\vert \varsigma^{*}_\zs{i}\vert]\right)
\\
&-
\sqrt{\frac{p^{*}_\zs{i}(1-p^{*}_\zs{i})}{u(1-u)}}
\,(u-p^{*}_\zs{i}) x\eta^{*}_\zs{i}(y,x)\,.
\end{align*}
\noindent
By Proposition \ref{Pr.sec:App.1-1} (see the appendix)
$
\max_\zs{1\le i\le N}\,
\sup_\zs{j\ge  1}\,
\EV^{*}\vert X_\zs{i,j}\vert^{\q}
<\infty
$
for any $\q>0$. Also,  under  $\Pb^{*}$ the random variables 
$(\eta^{*}_\zs{i}(X_\zs{i,j},X_\zs{i,j-1}))_\zs{1\le i\le N\,,\,j\ge 1}$ are i.i.d. $\Nc(0,1)$, so that
 condition $(\C^{*}_\zs{3}(\q))$ holds for any $\q>0$, which implies condition $(\C_\zs{3}(q))$ for any $\q>0$ (see Proposition \ref{Pr.sec:Mrk.1-1}).
Thus, it follows from Theorems \ref{Th.sec:Cnrsk.2} and  \ref{Th.sec:Cnrsk.2-*}
that the sequential detection-identification procedure $\delta_\beta^{opt}=(T_\beta^{opt}, d_\beta^{opt})$ defined in \eqref{seq-prs-opt}-\eqref{Topt}
is asymptotically optimal (as $\beta_{max} \to 0$) and robust in the pointwise and minimax senses for any $r>0$. 

\subsection{Monte Carlo} \label{MC_simulation} 

To get operating characteristics of the proposed detection-identification algorithm not only in the asymptotic case but also for reasonable probabilities of false alarm and misidentification,  
we perform Monte Carlo (MC) simulations for the modified Bernoulli model \eqref{ep-15-05-1}, \eqref{ep-15-05-2} with $X_\zs{i,n} = Y_\zs{i,n} /V_\zs{i}$.
The values of $Y_\zs{i,n}$ correspond to the number of susceptible at the $n$-th point in time for the 
$i$th population ($n \ge 0$, $i=1,\dots,N$) and the values of $V_\zs{i}$ to the number of susceptible at the initial moment, i.e. $Y_\zs{i,0} = V_\zs{i}$. 
In simulations, we set the initial value $Y_\zs{i,0} = V_\zs{i} = 0.5 (i + 1) \cdot 10^{4}$ for $1 \le i \le N$.
Without loss of generality we assume that the change occurs in the $N$th stream. Then \eqref{ep-15-05-1} for $1 \le i \le N - 1$ reduces to
$$
X_\zs{i,n}=(1-p^{*}_\zs{i})  X_\zs{i,n-1}+\sigma^{*}_\zs{i}\sqrt{\vert X_\zs{i,n-1}\vert}\,\xi_\zs{i,n}\,,\quad
 \sigma^{*}_\zs{i}=\sqrt{\frac{p_i^{*}(1-p_i^{*})}{V_\zs{i}}}\,,\quad X_\zs{i,0} = 1,
$$
and \eqref{ep-15-05-2} to
$$
X_\zs{N,n}=(1-\vartheta)  X_\zs{i,n-1}+\sigma_\zs{\vartheta}\sqrt{\vert X_\zs{N,n-1}\vert}\,\xi_\zs{N,n}\,,\quad
 \sigma_\zs{\vartheta}=\sqrt{\frac{\vartheta(1-\vartheta)}{V_\zs{N}}}\,,\quad X_\zs{N,0} = 1,
$$
where $\vartheta=\vartheta_\zs{n}=p^{*}_\zs{N}+(\overline{p}-p^{*}_\zs{N})\Chi_\zs{\{n>\nu\}}.$

In each MC run $m$, using formulas \eqref{g-1-16-05-2-0}--\eqref{ST-RdMatr-EP-33} and
\eqref{sec:Cnrsk.7}, we get a pair $\delta^{*, m}_\zs{A}=(T^{*, m}_\zs{A},d^{*, m}_\zs{A})$ --- the stopping time and
the number of the stream where the change is detected ($m = 1, \dots,M$,  $M$ is the total number of MC runs). 

The theoretic estimate of the expected detection delay for $i=N$ is given by the second asymptotic formula in \eqref{sec:SRAO2} with $r=1$, i.e., 
$$
\Rc_{N,\nu,\theta_N}\approx  \max_\zs{1\le j\le N} \frac{\vert \log\beta_\zs{j,N} \vert}{\iota_\zs{N,j}(\theta_N)}
\quad\mbox{for}\quad \theta_N=\overline{p} 
\,.
$$
Since calculation of $\iota_\zs{N,j}(\theta_N)$ analytically is difficult we evaluate it using MC simulations. To this end, we first estimate 
the conditional informations
$$
J_\zs{i}(\vartheta,x)=\int_\zs{\bbr}\,g_\zs{i}(\vartheta,y,x)\,f_\zs{\vartheta}(y\vert x)\d y
\quad\mbox{and}\quad
J^{*}_\zs{i}(\vartheta,x)=\int_\zs{\bbr}\,g_\zs{i}(\vartheta,y,x)\,f^{*}_\zs{i}(y\vert x)\d y
$$
and then we calculate the Kullback-Leibler divergences by MC as 
\begin{equation}
\label{Kulback-Appr-11-0}
\overline{I}_\zs{N}(\vartheta)=\frac{1}{K}\sum^{\nu+K}_\zs{n=\nu+1}\,J_\zs{N}(\vartheta,X_\zs{N,n})
 \quad\mbox{and}\quad
\overline{I}^{*}_\zs{i}(\vartheta)=\frac{1}{K}\sum^{K}_\zs{n=1}\,J^{*}_\zs{i}(\vartheta,X_\zs{i,n})
\quad\mbox{for}\quad 1\le i\le N-1
\,.
\end{equation}
By the law of large numbers for Markov chains  the MC estimates $\overline{I}_\zs{i}(\vartheta)$ and $\overline{I}^{*}_\zs{i}(\vartheta)$ converge to the true  values $\overline{J}_\zs{i}(\vartheta)$ 
and $\overline{J}^{*}_\zs{i}(\vartheta)$ defined in \eqref{sec:Mrk.6-erg}.  Then \eqref{Lower-bound-1} and \eqref{KL-Inf-1} reduce to
\begin{equation}
\label{Kulback-Appr-11}
\iota_\zs{N,j}(\theta_N)
=
\overline{I}_\zs{N}(\theta_N)
\Chi_\zs{\{j=N\}}
+
(\overline{I}_\zs{N}(\theta_N)-
\max_\zs{\vartheta\in \Theta_\zs{j}}
\overline{I}^{*}_\zs{j}(\vartheta))
\Chi_\zs{\{j\neq N\}}
\,.
\end{equation}

The MC estimate of the expected detection delay (in the $N$th stream) is calculated from the formula:
$$
\widehat{\Rc}=\frac{\sum^{M}_\zs{m=1} (T^{*,m}_\zs{A}-\nu) \Chi_\zs{(T^{*,m}_\zs{A} > \nu)} \Chi_\zs{\{d^{*,m}_\zs{A}=N\}}}{\sum^{M}_\zs{m=1} \Chi_\zs{\{T^{*,m}_\zs{A}>\nu\}}}\,.
$$
In particular, for $\nu=0$, which is used in simulations, it reduces to
\[
\widehat{\Rc}= \frac{1}{M} \sum^{M}_\zs{m=1} T^{*,m}_\zs{A} \Chi_\zs{\{d^{*,m}_\zs{A}=N\}}\,.
\]

The MC estimate of the false alarm probability ($\nu = +\infty$) is:
$$
\wh{\Pb}_\zs{N}
=
\max_\zs{1\le \ell\le \k^{*}-\m^{*}}
\frac{\sum^{M}_\zs{m=1} \Chi_\zs{\{\ell \le T^{*,m}_\zs{A}< \ell+\m^{*}\}}\, \Chi_\zs{\{d^{*,m}_\zs{A}=N\}}}{\sum^{M}_\zs{m=1} \Chi_\zs{\{T^{*,m}_\zs{A}\ge \ell\}}},
$$
and the MC estimates of the miss identification probabilities are:
$$
\check{\Pb}_\zs{j,N}
=
\max_\zs{\nu < \ell\le \nu + \k^{*}}
\,
\frac{\sum^{M}_\zs{m=1}\, \Chi_\zs{\{T^{*,m}_\zs{A}>\ell\}}\Chi_\zs{\{d^{*,m}_\zs{A}=j\}}}{\sum^{M}_\zs{m=1} \Chi_\zs{\{T^{*,m}_\zs{A}>\ell\}}}
\quad\mbox{for}\quad
1\le j\le N-1
\,.
$$

In simulations, we assume that the number of streams $N = 5$; the parameters of the observed process are $p^{*}_\zs{i}=1/(100+i), 1/(50+i)$; $q = \overline{p}/{p^{*}_\zs{N}} = 1.1, \, 1.15, \, 1.2$; 
for calculation of thresholds $A_\zs{i,j}$ we use \eqref{sec:Up-Bnd.1} and \eqref{sec:Cnrsk.7} with 
$\beta_\zs{i,j}=\frac{\varepsilon}{i+j}$, $\varepsilon = 0.3,\, 0.1,\, 0.01$ and $\check{\k} = 2,\, 1.55,\, 1.23$. We also assume that the change occurs from the very beginning, 
i.e., at the time $\nu=0$, in which case $\vartheta=\overline{p}$.

The results are shown in Table \ref{tab_Epidemic_M_C_1} and 
Table \ref{tab_Epidemic_M_C_2}. 
It is seen that the detection-identification algorithm has good performance. Even for small 
false alarm and miss identification probabilities the average
detection delay is small. Therefore,  we  recommend  using  this algorithm in practice 
for the detection and localization of epidemics. Also, the asymptotic approximations for the average detection delay are quite accurate and, therefore, can be used for the evaluation of the 
performance of the detection-identification procedure in practice.

\begin{table}[!h!]
\caption{Operating characteristics of the detection-identification procedure for $p^{*}_\zs{i}=\frac{1}{100+i}$ (MC simulations with $10^5$ runs).}
\begin{center}
\begin{tabular}{|c|c|c||c|c|c|c|c||c|c|c|}
\hline 
\textbf{$\varepsilon$} & \textbf{$\check{\k}$} & \textbf{$q$} & \textbf{$\check{\Pb}_\zs{1,N}$}& \textbf{$\check{\Pb}_\zs{2,N}$}& \textbf{$\check{\Pb}_\zs{3,N}$}& \textbf{$\check{\Pb}_\zs{4,N}$}& \textbf{$\wh{\Pb}_\zs{N}$}& \textbf{$\widehat{\Rc}$} & \textbf{$\Rc_{N,0,\theta_N}$} \\
\hline 
{0.3} &{2} & {1.1} &0.0024 &0.0027 &0.0011 &0.0007 &0.00088 & 6.46 & 5.17 \\
\hline 
{0.3} &{2} & {1.15} &0.0018 &0.0041 &0.0020 &0.0007 &0.00154 & 3.32 & 2.95 \\
\hline
{0.3} &{2} & {1.2} &0.0036 &0.0091 &0.0044 &0.0021 &0.00459 & 2.02 & 2.03 \\
\hline 
\hline 
{0.1} &{1.55} & {1.1} &0.0009 &0.0014 &0.0008 &0.0003 &0.00028 & 7.52 & 6.95 \\
\hline 
{0.1} &{1.55} & {1.15} &0.0004 &0.0023 &0.0013 &0.001 &0.0007 & 3.75 & 3.96 \\
\hline
{0.1} &{1.55} & {1.2} &0.0014 &0.0056 &0.0023 &0.0013 &0.0023 & 2.26 & 2.72 \\
\hline 
\hline 
{0.01} &{1.23} & {1.1} &0.00016 &0.00062 &0.00014 &0.00011 & $<10^{-5}$ & 9.96 & 10.50 \\
\hline 
{0.01} &{1.23} & {1.15} &0.0001 &0.0006 &0.0004 &0.0002 &0.00014 & 4.78 & 5.99 \\
\hline
{0.01} &{1.23} & {1.2} &0.0002 &0.0019 &0.0006 &0.0005 &0.0008 & 2.77 & 4.12 \\
\hline 
\end{tabular}
\label{tab_Epidemic_M_C_1}
\end{center}
\end{table}

\begin{table}[!h!]
\caption{Operating characteristics of the detection-identification procedure for $p^{*}_\zs{i}=\frac{1}{50+i}$ (MC simulations with $10^5$ runs).}
\begin{center}
\begin{tabular}{|c|c|c||c|c|c|c|c||c|c|c|}
\hline 
\textbf{$\varepsilon$} & \textbf{$\check{\k}$} & \textbf{$q$} & \textbf{$\check{\Pb}_\zs{1,N}$}& \textbf{$\check{\Pb}_\zs{2,N}$}& \textbf{$\check{\Pb}_\zs{3,N}$}& \textbf{$\check{\Pb}_\zs{4,N}$}& \textbf{$\wh{\Pb}_\zs{N}$}& \textbf{$\widehat{\Rc}$} & \textbf{$\Rc_{N,0,\theta_N}$} \\
\hline 
{0.3} &{2} & {1.1} &0.00076 &0.00062 &0.00033 &0.00016 &0.0004 & 3.92 & 3.89  \\
\hline 
{0.3} &{2} & {1.15} &0.00078 &0.00162 &0.0004 &0.0003 &0.0007 & 2.03 & 2.29 \\
\hline
{0.3} &{2} & {1.2} &0.0073 &0.0078 &0.0019 &0.0009 &0.0045 & 1.24 & 1.65  \\
\hline 
\hline 
{0.1} &{1.55} & {1.1} &0.0002 &0.00022 &0.00004 &0.00006 &0.00014 & 4.50 & 5.23 \\
\hline 
{0.1} &{1.55} & {1.15} &0.00018 &0.00082 &0.00026 &0.00016 &0.00044 & 2.26 & 3.08  \\
\hline
{0.1} &{1.55} & {1.2} &0.00342 &0.00447 &0.00142 &0.00049 &0.0017 & 1.33 & 2.21 \\
\hline 
\hline 
{0.01} &{1.23} & {1.1} &0.00004 &0.00004 &0.00002 &0.00002 &0.00002 & 5.74 & 7.91 \\
\hline 
{0.01} &{1.23} & {1.15} &0.00002 &0.00022 &0.00008 &0.00001 &0.00008 & 2.78 & 4.65 \\
\hline
{0.01} &{1.23} & {1.2} &0.00084 &0.00145 &0.00056 &0.00019 &0.0003 & 1.57 & 3.35 \\
\hline 
\end{tabular}
\label{tab_Epidemic_M_C_2}
\end{center}
\end{table}

\subsection{Detection of COVID-19 in Italy}\label{ssec:COVID}

In Subsection \ref{ssec:AO}, we applied the proposed sequential detection-identification algorithm to epidemic models and showed it to be asymptotically optimal when the 
probabilities of wrong identification and false alarm are small. In this subsection, we demonstrate that the proposed detection-identification procedure can be effectively applied for the localization of 
COVID-19, i.e., for the detection of the epidemic anomalies and identification of the affected region. Consider the case of Italy. 

Let $H_\zs{i,n}$ 
be the number of hospitalized people at the $n$-th moment for the $i$th region. (Since the shortage of hospital beds presented a major challenge in Italy during the first wave of COVID-19, 
we focus on hospitalizations. However, this model also applies to other kinds of observations, e.g., number of infected people, number of visits to the doctor \cite{BaronChoudharyYu2013}.) Then, 
$Y_\zs{i,n} = V_\zs{i} - H_\zs{i,n},$ where $V_\zs{i}$ is the total number of hospital beds, i.e., $Y_\zs{i,n}$ 
is potentially free beds for new hospitalizations at the $n$-th moment for the $i$th region. Then the observation, 
as in \eqref{MC_simulation},  will be $X_\zs{i,n} = Y_\zs{i,n} / V_\zs{i}$.

\begin{figure}[!h!] \centering
\includegraphics[scale=0.65]{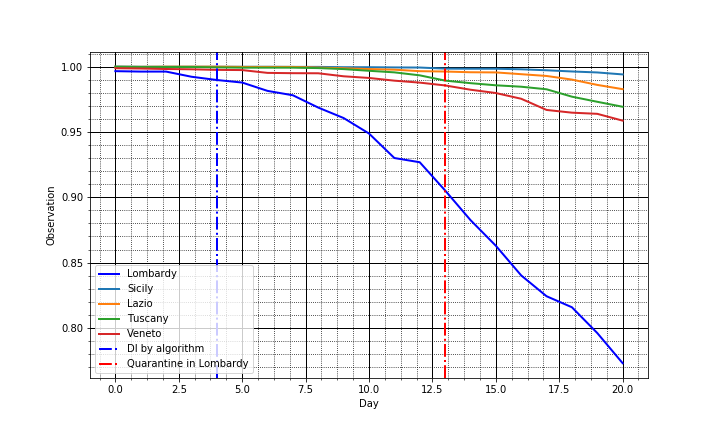}
\caption{Detection and identification of the region with an outbreak of the COVID-19 epidemic in Italy: the proposed algorithm vs. the imposition of quarantine protocols in Lombardy. The time of detection and
identification by our algorithm is shown by blue vertical line and the time of the official imposition of quarantine protocols by red vertical line. }
\label{Paper_Epidemic_Italy}
\end{figure}

We use the data provided by {\em Sito del Dipartimento della Protezione Civile - Emergenza Coronavirus: la risposta nazionale} (the Italian Department of Civil Protection). This data includes 
information on hospitalizations by region each day. We consider five Italian regions: Sicily, Lazio, Tuscany, Venice, and Lombardy. 
 We use the proposed detection-identification algorithm to detect the presence of COVID-19 in a given region. 
Fig. \ref{Paper_Epidemic_Italy} shows raw observations for five different
regions, detection and identification of a region with a COVID outbreak in Italy by the proposed algorithm (blue vertical line), and the official introduction of a regional quarantine (red vertical line) in Lombardy.  

It is known that Lombardy became the epicenter of the spread of COVID not only in Italy but throughout Europe. 
According to Fig. \ref{Paper_Epidemic_Italy}, the proposed algorithm detected COVID in Italy 9 days prior to the imposition of quarantine protocols in Lombardy (February 28, 2020 vs. March 8, 2020). The proposed detection-identification algorithm could therefore be a useful tool for researchers and public health aurhorities in detecting and localizing epidemics.

\section{Conclusion}\label{sec:Conclusion}

1. In this paper, we ignore the possible indifference zone $\Theta_{\rm ind}$ of parameter values where the probabilities of false alarms and misidentification are too close to be reasonably distinguishable. 
In the indifference zone, the constraints on the erroneous decisions are not imposed, but still, the expected detection delays (or more generally moments of delay) have to be minimized for all 
possible parameter values, including those in the indifference zone. The modification of the
proposed procedure to take into account an indifference zone, if needed, is straightforward. For the sake of brevity, the details are omitted.

2. As in the recent paper by Tartakovsky~\cite{TartakovskyIEEEIT2021}, we focus on the multistream changepoint model \eqref{des-1}--\eqref{dens-2}. It is worth noting that the same 
results hold in the single-stream detection-isolation problem when the 
observations $\{X_n\}_{n \ge 1}$  represent either a scalar process or a vector process but 
all components of this process change at time $\nu$. Specifically, in change detection and isolation, the post-change hypothesis $H_{\nu,i,\theta_i}, \theta\in \Theta_i$ corresponding to the $i$th type of change
usually involves unknown parameters $\theta_i$ and, therefore, is composite. Under the hypothesis $H_{\nu,i,\theta_i}, \theta\in \Theta_i$ the post-change conditional density function is
$f_{i,\theta_i,n}(X_n \vert \Xb_1^{n-1})$, $n >\nu$, while the pre-change density is $f^*_{n}(X_n \vert \Xb^{n-1})$, $n \le \nu$, where $\Xb^t=(X_1,\dots,X_t)$.
Hence, introducing parametric families of densities  $\{f_{i,\theta_i,n}(X_n|\Xb^{n-1}), \theta_i\in\Theta_i\}$
and for  $i=1,\dots,N$ and $\Theta_i\subset\Theta$ considering the model\footnote{Often, $\theta_i\equiv \theta$ does not depend on $i$ in practice.}  
\begin{equation} \label{density}
p (\Xb^n \vert H_{\nu,i,\theta_i}) =
\begin{cases}
 \q^{*}(\Xb^n) & \quad \text{for}~~ \nu \ge n\,;
\\ 
\q^{*}(\Xb^\nu) \, \prod_{l=\nu+1}^{n} f_\zs{i,\theta_\zs{i},l}(X_\zs{l}|\Xb^{l-1})&  \quad \text{for}~~ \nu < n
\end{cases} ,
\end{equation}
where $p (\Xb^n \vert H_{\nu,i,\theta_i})$ stands for the joint density of the first $n$ observations $\Xb^n$ conditioned on the hypothesis $H_{\nu,i,\theta_i}$ and
$\q^{*}(\Xb^s) = \prod_\zs{l=1}^s f^{*}_\zs{l}(X_\zs{l}| \Xb^{l-1})$ for $s \ge 1$,
we arrive at the single-stream model that has all features of the previous multistream model \eqref{des-1}--\eqref{dens-2}. In fact, setting $X_n=(X_{1,n},\dots,X_{N,n})$, where the components of this vector
are mutually independent and assuming that the change may occur only in a single component, we obtain 
\begin{equation*}
p(\Xb^n \vert H_{\nu,i,\theta_i})  =
\begin{cases}
  \prod_{i=1}^N \q^{*}_i(\Xb_i^n) & \quad \text{for}~~ \nu \ge n\,;
\\
 \q^{*}_i(\Xb_i^\nu) \, \prod_{l=\nu+1}^{n}  f_{i,\theta_i,l}(X_{i,l}|\Xb_i^{l-1}) \prod_{1\le j\neq i \le N}   \prod_{l=1}^n f^*_j(X_{j,l})|\Xb_j^{l-1})  & \quad \text{for}~~ \nu < n
\end{cases} ,
\end{equation*}
where $\Xb_i^n=(X_{1,i},\dots,X_{i,n})$, $\Xb^n=(\Xb_1^n,\dots,\Xb_N^n)$, and  $\q_i^{*}(\Xb^s) = \prod_\zs{l=1}^s f^{*}_\zs{i,l}(X_\zs{l}| \Xb_i^{l-1})$. Obviously, this joint density is the same as the one
in \eqref{jointforall}, so that in the case of mutually independent streams the multistream model defined in \eqref{des-1}--\eqref{dens-2} is a particular case of the model \eqref{density}.

3.  All previous results can be generalized for the case when the change points are different for different streams, i.e., when $\nu=\nu_i$. 

4. For independent observations as well as for a variety of Markov and certain hidden Markov models (see, e.g., Subsections~\ref{ssec:Ex1} and \eqref{ssec:Ex2} and Section~\ref{sec:Epid}), 
the decision statistics $<\U_\zs{n}>_\zs{i,j}$ defined in \eqref{ST-RdMatr} can be computed relatively easily, in which case 
implementation of the proposed change detection-identification procedure is not
an issue. In general, however, the computational complexity of rule $\delta_A$ may be high. To avoid computational difficulties 
rule $\delta_A$ can be modified into a window-limited version where the summation in the decision statistics over potential change points $\nu$ is restricted to the sliding window of specific fixed size 
$\ell=\ell_\beta$, which is a function of the error probabilities constraints $\beta$.  Following guidelines of \cite{TartakovskyIEEEIT2021} (Ch 3, Sec 3.10)
where asymptotic optimality of mixture window-limited rules was established in the single-stream case, it can be shown that the window-limited version of the multihypothesis detection-identification procedure
also has asymptotic optimality properties as long as the size of the  window $\ell_\beta$ goes to infinity as $\beta_{max} \to 0$ at such a rate that $\ell_\beta /|\log\beta_{max}| \to \infty$.

\section*{Acknowledgements}

The work of S.M. Pergamenchtchikov  was  partially supported by  the RSF  grant 20-61-47043  (National Research Tomsk State University).  
The work of A.G. Tartakovsky  was supported in part by  the Russian Science Foundation  Grant 18-19-00452 at the Moscow Institute of Physics and Technology (Space Informatics Laboratory).

\appendix



\renewcommand{\thetheorem}{A.\arabic{theorem}}
\setcounter{theorem}{0}

\renewcommand{\theproposition}{A.\arabic{proposition}}
\setcounter{proposition}{0}

\section{Moment properties of the epidemic models }\label{subsec:Mprts}

\begin{proposition}
\label{Pr.sec:App.1-1}
For any integer $m\ge 1$ and for any compact set $\Theta_\zs{i}\subseteq \Theta$ the process $X_{i,n}$ defined in \eqref{ep-15-05-1} and \eqref{ep-15-05-2} 
has the following moment properties
\begin{equation}
\label{Upper-Bnd-16-05-1}
\sup_\zs{x\in\bbr}
\sup_\zs{n\ge 0}
\max_\zs{1\le i\le N}
\frac{\EV^{*}_\zs{i,x}\,(X^{*}_\zs{i,n})^{2m}}{1+x^{2m}}
<\infty
\quad\mbox{and}\quad
\sup_\zs{x\in\bbr}
\sup_\zs{n\ge 0}
\max_\zs{1\le i\le N}
\sup_\zs{\vartheta\in\Theta_\zs{i}}
\frac{\EV^{\vartheta}_\zs{i,x}\,(X_\zs{i,n})^{2m}}{1+x^{2m}}
<\infty
\,.
\end{equation}

\end{proposition}

\proof
We prove only the second inequality in \eqref{Upper-Bnd-16-05-1} since the proof of the first one is essentially similar.
To this end, we first show that for any $x\in\bbr$
\begin{equation}
\label{Upper-Bnd-16-05-1-2}
\sup_\zs{n\ge 0}
\max_\zs{1\le i\le N}
\sup_\zs{\vartheta\in\Theta_\zs{i}}
\EV^{\vartheta}_\zs{i,x}\,[X^{2}_\zs{i,n}]
\le x^{2}+1
\,.
\end{equation}
It is easily seen that for the model \eqref{ep-15-05-2} we have 
$$
\EV^{\vartheta}_\zs{i,x}\,[X^{2}_\zs{i,n}]\le (1-\vartheta)^{2}\EV^{\vartheta}_\zs{i,x}\,[X^{2}_\zs{i,n-1}]
+\sigma^{2}_\zs{\vartheta} \EV^{\vartheta}_\zs{i,x}\,\vert X_\zs{i,n-1}\vert
\,.
$$
For the sake of brevity write $y_\zs{n}=\EV^{\vartheta}_\zs{i,x}\,[X^{2}_\zs{i,n}]$. Taking into account that $x\le 1+x^{2}$, we obtain 
$$
y_\zs{n}\le 
(1-\vartheta)^{2}y_\zs{n-1}+\sigma^{2}_\zs{\vartheta}\sqrt{y_\zs{n-1}}
\le 
[(1-\vartheta)^{2}+\sigma^{2}_\zs{\vartheta}]y_\zs{n-1}
+\sigma^{2}_\zs{\vartheta}
= (1-\vartheta) y_\zs{n-1}
+\sigma^{2}_\zs{\vartheta}
\,.
$$
Let now $\upsilon_\zs{n}=y_\zs{n}-(1-\vartheta)^{2}y_\zs{n-1}$. Clearly $\upsilon_\zs{n}\le \sigma^{2}_\zs{\vartheta}$
and
$$
y_\zs{n}=x^{2}\,(1-\vartheta)^{n}+\sum^{n}_\zs{j=1}\vartheta^{n-j}\,\upsilon_\zs{j}
\le x^{2}+\sigma^{2}_\zs{\vartheta}\sum_\zs{j\ge 0}\,(1-\vartheta)^{j}
= 
x^{2}+\frac{\sigma^{2}_\zs{\vartheta}}{\vartheta}
=x^{2}+1-\vartheta\,.
$$
This implies inequality \eqref{Upper-Bnd-16-05-1-2} and, therefore, the second inequality in  \eqref{Upper-Bnd-16-05-1} for $m=1$.

For an arbitrary $m\ge 1$ this inequality can be proved by induction as follows. Assume that the second inequality in 
\eqref{Upper-Bnd-16-05-1} is true for $m-1$ and 
$m\ge 2$, i.e., there exists a constant $C_\zs{m}\ge 1$ such that for any
$x\in\bbr$, $k\ge 1$, $1\le i\le N$ and $\vartheta\in\Theta_\zs{i}$
\begin{equation}
\label{sec:App-Ind-1}
\EV^{\vartheta}_\zs{i,x}\,[(X_\zs{i,k})^{2(m-1)}]
\le 
C_\zs{m}(1+x^{2(m-1)})
\,.
\end{equation}
To show that it holds for $m$, using the initial condition $X_\zs{i,0}=x$,  we represent the process \eqref{ep-15-05-2} as
$$
X_\zs{i,k}=(1-\vartheta)^{k}\,x+ \sigma_\zs{\vartheta} \sum^{k}_\zs{j=1}(1-\vartheta)^{k-j}\,\sqrt{\vert X_\zs{i,j-1}\vert}\xi_\zs{i,j}
\,.
$$
By the H\"older inequality, 
\begin{align*}
\EV^{\vartheta}_\zs{i,x}\,[X^{2m}_\zs{i,k}]&\le 2^{2m-1}\left( 
x^{2m}
+\frac{\sigma^{2m}_\zs{\vartheta}}{\vartheta^{2m-1}}
\sum^{k}_\zs{j=1}\,
(1-\vartheta)^{k-j} \EV^{\vartheta}_\zs{i,x}[\vert X_\zs{i,k}\vert^{m}\xi^{2m}_\zs{i,j}]
\right)\\[3mm]
&\le 
2^{2m-1}\left( 
x^{2m}
+\frac{(2m-1)!!\sigma^{2m}_\zs{\vartheta}}{\vartheta^{2m}}
 \sup_\zs{k\ge 0}\EV^{\vartheta}_\zs{i,x}\vert X_\zs{i,k}\vert^{m}
\right)\,.
\end{align*}
 \noindent
Now, using the induction assumption \eqref{sec:App-Ind-1}
and that $|x+y|^{\alpha}\le |x|^{\alpha}+|y|^{\alpha}$ for $0<\alpha\le 1$,
 we obtain
$$
\EV^{\vartheta}_\zs{i,x}\vert X_\zs{i,k}\vert^{m}\le 
\left(
\EV^{\vartheta}_\zs{i,x}\vert X_\zs{i,k}\vert^{2(m-1)}
\right)^{m/(2m-2)}
\le C^{m/(2m-2)}_\zs{m}
\left(
1+
x^{2(m-1)}
\right)^{m/(2m-2)}
\le 
C_\zs{m}
\left(
1+
x^{2m}
\right)
\,.
$$
This implies the second inequality in\eqref{Upper-Bnd-16-05-1}, completing the proof.
\endproof





\end{document}